\documentclass[11pt,english]{article}
\usepackage[left=1in,right=1in,top=1.in,bottom=1.in]{geometry}
\usepackage{authblk}
\usepackage{graphicx,graphics}
\usepackage{amssymb}
\usepackage{amsbsy}
\usepackage{stmaryrd} 
\usepackage{amsmath}
\usepackage{multirow}
\usepackage{graphicx}
\usepackage{url}  
\usepackage[toc,page]{appendix}

\usepackage{ragged2e}
\usepackage{dblfloatfix}
\usepackage{booktabs}       

\DeclareMathOperator*{\argmin}{argmin}
\DeclareMathOperator*{\argmax}{argmax}
\DeclareMathOperator*{\Supp}{Supp}

\newtheorem{theorem}{Theorem}
\newtheorem{defn}{Definition}
\newtheorem{lemma}{Lemma}
\newtheorem{cor}{Corollary}
\newtheorem{asu}{Assumption}
\newenvironment{Proof}{\paragraph{Proof:}}{\hfill$\square$}

\usepackage{natbib}

\newcommand{\sbt}{\mathrm{s.t. }}
\newcommand{\sign}{\text{sign}}
\newcommand{\subg}{\text{subG}}
\newcommand{\B}{\boldsymbol}

\newenvironment{myarray}[2][1]
{\def\arraystretch{#1}\array{#2}}
{\endarray}

\usepackage[toc,page]{appendix}
\usepackage{times}

\begin{document}

\author{Antoine Dedieu \thanks{Antoine Dedieu’s research was partially supported by the Office of Naval Research: N000141512342 and by Vicarious. Email: \texttt{ adedieu@mit.edu, antoine@vicarious.com}.}}
\affil{Massachusetts Institute of Technology, Vicarious}

\date{\today}

\title{Error bounds for sparse classifiers in high-dimensions}
\maketitle

\begin{abstract}
	We prove an L2 recovery bound for a family of sparse estimators defined as minimizers of some empirical loss functions -- which include hinge loss and logistic loss. More precisely, we achieve an upper-bound for coefficients estimation scaling as $(k^*/n)\log(p/k^*)$: $n\times p$ is the size of the design matrix and $k^*$ the dimension of the theoretical loss minimizer. This is done under standard assumptions, for which we derive stronger versions of a cone condition and a restricted strong convexity. Our bound holds with high probability and in expectation and applies to an L1-regularized estimator and to a recently introduced Slope estimator, which we generalize for classification problems.  Slope presents the advantage of adapting to unknown sparsity. Thus, we propose a tractable proximal algorithm to compute it and assess its empirical performance. Our results match the best existing bounds for classification and regression problems. 
\end{abstract}

\section{Introduction}
Motivated by the increasing availability of very large-scale datasets, high-dimensional statistics has focused on analyzing the performance of sparse estimators.  An estimator is said to be sparse if the response of an observation is given by a small number of coefficients: sparsity delivers better interpretability and often leads to computational efficiency. Statistical performance and L2 consistency for high-dimensional linear regression have been widely studied. For two polynomial-time sparse estimators, a Lasso  \citep{tibshirani1996regression} and a Dantzig selector \citep{candes2007dantzig},  \citet{lasso-dantzig}  proved a $(k^*/n)\log(p)$ rate for the L2 estimation of the coefficients: $n\times p$ is the dimension of the input matrix and $k^*$ the degree of sparsity of the vector used to generate the model. The optimality of this bound is essential for a theoretical understanding of the method performance.  \citet{candes-sparse-estimation} and  \citet{raskutti_wainwright} proved a  $(k^*/n)\log(p/k^*)$ lower bound for estimating the L2 norm of a sparse vector, regardless of the input matrix and estimation procedure. This optimal minimax rate is known to be achieved by a sparse but theoretically intractable BIC estimator  \citep{bic-tsybakov} which considers an L0 regularization. The BIC estimator adapts to unknown sparsity: the degree $k^*$ does not have to be specified. Recently, \citet{slope} reached this optimal minimax bound for a Lasso estimator with knowledge of the sparsity $k^*$. They also proved that a recently introduced and polynomial-time Slope estimator \citep{slope-introduction} achieves this optimal rate while adapting to unknown sparsity.

\medskip

Little work has been done on deriving (theoretical) upper bounds for the estimation error on high-dimensional classification problems: the literature has essentially focused on analysis of convergence \citep{tarigan, vssvm}. Recently, \citet{L1-SVM} proved a $(k^*/n) \log(p)$ upper-bound for L2 coefficients estimation of a L1-regularized Support Vector Machines (SVM): $k^*$ is now the sparsity of the theoretical minimizer to estimate. They recovered the rate proposed by \citet{vdg_linear_models}, which considered a weighted L1 norm for linear models. \citet{Wainwright-logreg} obtained a similar bound for a L1-regularized Logistic Regression estimator in a binary Ising graph. Their frameworks and bounds are similar to the model proposed by \citet{quantile-reg} for L1-regularized Quantile Regression; this inspired us to include this problem in our framework. However, this rate of $(k^*/n)\log(p)$ is not the best known for a classification estimator: \citet{one-bit} proved a $k^*\log(p/k^*)$ error bound for estimating a single vector through sparse models -- including 1-bit compressed sensing and Logistic Regression -- over a bounded set of vectors. Contrary to this work, our approach does not assume a generative vector and applies to a larger class of problems (SVM, Quantile Regression) and regularizations (Slope). In addition, our framework share similarity with Section 4.4. of \citet{M-estimators}: the authors consider some sub-gaussian tails assumptions and restricted eigenvalue conditions to derive a restricted strong convexity condition similar to our Theorem \ref{restricted-strong-convexity}. However, their results only apply to generalized linear models, and are weaker since the parameter $\tau(k)$ proposed in the tolerance function of the restricted strong convexity condition is higher than ours. Finally, \citet{pierre2017estimation} studied a similar class of loss functions and regularization that the ones proposed herein. However, their proof technique is quite different than ours, leading to an estimation error rate of the order of $k^*\log(p)/n$, which is higher than the one we derive. The authors do not discuss any computational algorithms for the Slope estimator, which we do. 

\medskip

\noindent{\bf{What this paper is about:} } In this paper, we propose a theoretical framework to analyze the properties of a general class of sparse estimators for classification problems -- including SVM and Logistic Regression -- with different regularization schemes. Our approach draws inspiration from the least squares regression case and illustrates the distinction between regression and classification studies. Our main results are first presented for a family of L1-regularized estimators. We achieve a $(k^*/n)\log(p/k^*)$ upper-bound for coefficients estimation, which holds with high probability and in expectation. In addition, we introduce a version of the Slope estimator for classification problems: we propose a proximal algorithm to compute the solution, and we prove that a tractable Slope estimator achieves a similar upper-bound while adapting to unknown sparsity. To the best of our knowledge, it is the first time any of these bounds is reached for the estimators considered.

\medskip

The rest of this paper is organized as follows. Section \ref{sec: framework}  introduces and discusses common assumptions in the literature, and builds our framework of study in the case of L1-regularized estimators. Section \ref{sec: error-bound} proves two essential results and derive our upper-bounds in Theorem \ref{main-results} and Corollary \ref{main-corollary}. Finally, Section \ref{sec: slope} defines and computes the Slope estimator for our class of problems and discusses its statistical performance.

\section{General assumptions with an L1 regularization}\label{sec: framework} 

We consider a set of training data $\left\{ (\mathbf{x_i},y_i) \right\}_{i=1}^n$, $( \mathbf{x_i},y_i) \in \mathbb{R}^p\times \mathcal{Y}$ from an unknown distribution $\mathbb{P}(\mathbf{X}, \mathbf{y} )$. We note our loss $f$ and define the theoretical loss $\mathcal{L}(\B{\beta}) = \mathbb{E} \left( f \left( \langle \mathbf{x},  \B{\beta} \rangle ;  y \right)  \right)$. We consider a theoretical minimizer $ \B{\beta}^*$:
\begin{equation} \label{def-beta0}
\B{\beta}^* \in \argmin \limits_{ \B{\beta} \in \mathbb{R}^{p}} \left\{\mathbb{E} \left( f \left( \langle \mathbf{x},  \B{\beta} \rangle ;  y \right) \right)  \right\}.
\end{equation}
In the rest of this section, we denote by $k = \| \B{\beta}^*\|_0$ the number of non-zeros of the theoretical minimizer and $R=\| \B{\beta}^*\|_1$ its L1 norm. We assume $R\ge 1$. We study the L1-regularized L1-constrained problem defined as:
\begin{equation} \label{learning}
\min \limits_{ \B{\beta} \in \mathbb{R}^{p}:\  \| \B{\beta}  \|_1 \le 2R } \;\; \frac{1}{n}  \sum_{i=1}^n f \left( \langle \mathbf{x_i},  \B{\beta} \rangle ;  y_i \right) + \lambda \| \B{\beta} \|_1.
\end{equation}
We consider an empirical minimizer $\hat{\B{\beta}}$, solution of Problem \eqref{learning}. The constraint $2R$ in Problem \eqref{learning} is somewhat arbitrary: it enforces the empirical minimizer to be close enough to the theoretical minimizer $\B{\beta}^*$: $\| \hat{\B{\beta}}- \B{\beta}^* \|_1 \le 3R$. The L1 regularization in Lagrangian form is known to induce sparsity in the coefficients of  $\hat{\B{\beta}}$. Note that Problem \eqref{learning} is fully tractable.

\medskip

For a given $\lambda$, we fix a solution $\hat{\B{\beta}}(\lambda,R) $ of Problem \eqref{learning} -- $R$ is fixed throughout the paper. Our main result is an error bound -- achieved for a certain $\lambda$ -- for the L2 norm of the difference between the empirical and theoretical minimizers $\| \hat{\B{\beta}}(\lambda,R) - \B{\beta}^*\|_2$. When no confusion can be made, we drop the dependence upon the parameters $\lambda,R$. Our bound is reached under standard assumptions in the literature. In particular, it is similar to those proposed by \cite{L1-SVM}, \cite{Wainwright-logreg}, \cite{quantile-reg}. The rest of this section presents our framework of study.

\subsection{Lipschitz loss function}
Our first assumption concerns the Lipschitz-continuity of the loss $f$. 
\begin{asu} \label{asu1}
	The loss $f(. , y)$ is non-negative, convex and Lipschitz continuous with constant $L$, that is, $| f(t_1, y) - f(t_2, y) | \le L | t_1 -t_2 |, \ \forall t_1, t_2$.
	In addition, there exists $\partial f(.,y)$ such that $f(t_2, y) - f(t_1, y) \ge \partial f(t_1,y) (t_2 -t_1), \ \forall t_1, t_2$. 
\end{asu}
$\partial f(.,y)$ is said to be a sub-gradient of the loss: if $f(. , y)$ is differentiable, we simply consider its gradient. It trivially holds $\| \partial f(.,y) \|_{\infty} \le L, \ \forall y$. We list three main examples that fall into this framework.

\paragraph{Example 1: Support Vectors Machines}
We assume $\mathcal{Y} = \left\{ -1, 1\right\} $ and consider the L1-regularized L1-constrained Support Vector Machines (SVM) problem. It learns a classification rule of the data of the form $\sign( \langle \mathbf{x},  \B{\beta} \rangle )$ by solving the problem:
\begin{equation} \label{SVM}
\min \limits_{ \B{\beta} \in \mathbb{R}^{p}:\  \| \B{\beta}  \|_1 \le 2R } \;\; \frac{1}{n}  \sum_{i=1}^n \left( 1 - y_i \langle \mathbf{x}_i,  \B{\beta} \rangle  \right)_+ + \lambda \| \B{\beta} \|_1.
\end{equation}
The hinge loss $f \left( \langle \mathbf{x},  \B{\beta} \rangle ;  y \right) = \max(0, 1 - y \langle \mathbf{x},  \B{\beta} \rangle )$ admits as a subgradient $\partial f(.,y)=\mathbf{1}(1-y. \ge 0)y.$ and satisfies Assumption \ref{asu1} for $L=1$. 

\paragraph{Example 2: Logistic Regression}
Here, we still have $\mathcal{Y} = \left\{ -1, 1\right\} $ and we consider the additional assumption $\log\left( \mathbb{P}(y_i =1 | \mathbf{X} =  \mathbf{x_i}) \right) - \log\left( \mathbb{P}(y_i = -1 | \mathbf{X} =  \mathbf{x_i}) \right)= \langle \mathbf{x}_i,  \B{\beta} \rangle,  \ \forall i $.  The  L1-regularized L1-constrained Logistic Regression estimator is a solution of the problem:
\begin{equation} \label{logreg-bis}
\min \limits_{ \B{\beta} \in \mathbb{R}^{p}:\  \| \B{\beta}  \|_1 \le 2R } \frac{1}{n}  \sum_{i=1}^n \log \left( 1 + \exp ( - y_i\langle \mathbf{x}_i,  \B{\beta} \rangle  ) \right) + \lambda \| \B{\beta} \|_1.
\end{equation}
The logistic loss $f \left( \langle \mathbf{x},  \B{\beta} \rangle ;  y \right) = \log(1 + \exp ( - y \langle \mathbf{x},  \B{\beta} \rangle ))$ has a derivative with repect to its first variable  $| \partial_t f(t,y) | = \left|1/ \left(1+ e^{yt} \right) \right| \le 1$, hence it satisfies Assumption \ref{asu1} for $L=1$.

\paragraph{Example 3: Quantile Regression}
We now consider a class of parametric quantile estimation problems. Following \cite{advance-quantile-reg}, we assume that for  $\theta \in (0,1)$ the conditional quantile of $y$ given $\mathbf{X}$ is given by $Q_{\theta} (y | \mathbf{X} =  \mathbf{x}) = \langle \mathbf{x},  \B{\beta}_{\theta} \rangle$, where the model is of the form $y =  \langle \mathbf{x},  \B{\beta}_{\theta} \rangle + u_{\theta}$, and $u_{\theta}$ is unkown. The L1-regularized L1-constrained $\theta$-Quantile Regression estimator is  defined as a solution of:
\begin{equation} \label{quantile-reg}
\min \limits_{ \B{\beta} \in \mathbb{R}^{p}:\  \| \B{\beta}  \|_1 \le 2R  } \frac{1}{n}  \sum_{i=1}^n \rho_{\theta} \left( y_i - \langle \mathbf{x}_i,  \B{\beta} \rangle  ) \right) + \lambda \| \B{\beta} \|_1,
\end{equation}
where $\rho_{\theta} (t) = (\theta - \mathbf{1}(t\le0 ) )t$ is the quantile regression loss. $\rho_{\theta}$ satisfies Assumption \ref{asu1} for $L = \max(1 - \theta, \theta).$ Note that the hinge loss is a simple translation of the quantile regression loss for $\theta = 0$.

\subsection{Differentiability of the theoretical loss}\label{sec: TheoreticalMin}

The following assumption ensures the unicity of $ \B{\beta}^*$ and the twice differentiability of the theoretical loss $\mathcal{L}$. Equation \eqref{gradient-relation} is equivalent to saying that the gradient of the theoretical loss is equal to the theoretical sub-gradient of the loss -- defined in Assumption \ref{asu1}.
\begin{asu} \label{asu2}
	The theoretical minimizer is unique. In addition, the theoretical loss is twice-differentiable: we denote its gradient $\nabla \mathcal{L}(\B{\beta}) $ and its Hessian matrix $\nabla^2 \mathcal{L}(\B{\beta}).$ We also assume:
	\begin{equation}\label{gradient-relation}
	\nabla \mathcal{L}(.) = \mathbb{E}\left( \partial f \left( \langle \mathbf{x},  .  \rangle ;  y \right) \mathbf{x} \right).
	\end{equation}
\end{asu}

\paragraph{Support Vectors Machines:}
\cite{lemma2} studied specific conditions under which Assumption \ref{asu2} holds for SVM. In particular, if $f$ and $g$ denote  the respective conditional densities of  $\mathbf{X}$ given $y=1$ and $y=-1$; they proved that if the densities $f$ and $g$ are continuous with common support $\mathcal{S} \subset \mathbb{R}^p$ and have finite second moments, then the gradient $\nabla \mathcal{L}(\B{\beta} ) = \mathbb{E}\left(  \mathbf{1} \left(  1 -  y \langle \mathbf{x},  \B{\beta}   \rangle  \ge 0\right) y \mathbf{x} \right)$ and the Hessian matrix $\nabla^2 \mathcal{L}(\B{\beta} ) = \mathbb{E}\left(  \delta \left(  1 -  y \langle \mathbf{x},  \B{\beta}   \rangle \right) y \mathbf{x} \right)$ ( $\delta(.)$ is the Dirac function ) are defined and continuous.

\paragraph{Logistic and Quantile Regression:} 
The regularity of $\nabla \mathcal{L}$ and $\nabla^2 \mathcal{L}$ are trivial for the logistic regression loss. Equation \eqref{gradient-relation} holds as the sub-gradient is simply the gradient of the loss.  For the quantile regression loss, a study similar to the case of the hinge loss -- using Assumption D.1 by \cite{quantile-reg} -- can be applied to obtain Assumption \ref{asu2}.

\subsection{Sub-Gaussian columns}\label{sec: asu-SG}
We denote  $\mathbb{X}$ the design matrix, with rows $\mathbf{x}_1, \ldots, \mathbf{x}_n$. The following assumption guarantees that some random variables of the columns $ (\mathbf{X}_1, \ldots, \mathbf{X}_p)$ of  $\mathbb{X}$ have their tails bounded by a sub-Gaussian random variable with variance proportional to $n$. We first recall the definition of a sub-Gaussian random variable~\citep{lecture-notes}:
\begin{defn} \label{def-asu3}
	A random variable $Z$ is said to be sub-Gaussian with variance $\sigma^2>0$ if $\mathbb{E}(Z) = 0$ and $\mathbb{P}\left( |Z| > t \right) \le 2 \exp \left( - \frac{t^2}{2\sigma^2} \right), \  \forall t>0$.
\end{defn}
A sub-Gaussian variable will be noted  $Z \sim \subg(\sigma^2)$. We would like here to notice another important aspect of our contribution. Our next Theorem \ref{sec: error-bound} derives a cone condition, a necesary step to prove our main results. Our approach draws inspiration from the regression case with Gaussian noise. However, it relies on a new study of sub-Gaussian random variables -- such analysis is not needed in the regression case. Our results are derived under the following Assumption \ref{asuSG}:
\begin{asu} \label{asuSG}
	There exists $M>0$ such that with the notations of Assumption \ref{asu1}:
	\begin{equation}
	\sum_{i=1}^n  \partial f\left( \langle \mathbf{x}_i,  \B{\beta}^*  \rangle, y_i  \right)  x_{ij}\sim \subg(n L^2 M^2), \ \forall j.
	\end{equation}
\end{asu}
$\B{\beta}^*$ minimizes the theoretical loss. Thus, from Assumption \ref{asu2},  $\mathbb{E} \left[  \partial f\left( \langle \mathbf{x}_i,  \B{\beta}^*  \rangle, y_i  \right)  x_{ij}\ \right] =0, \forall i,j$. The next lemma gives more insight about Assumption \ref{asuSG}. The proof is presented in Appendix \ref{sec: appendix_lemma}.
\begin{lemma} \label{lemma-asu3}
	If the rows of the design matrix are independent and if all the entries $\partial f\left( \langle \mathbf{x}_i,  \B{\beta}^*  \rangle; y_i  \right) x_{ij}, \ \forall i,j$ are sub-Gaussian with variance $L^2M^2$, then $\sum_{i=1}^n  \partial f\left( \langle \mathbf{x}_i,  \B{\beta}^*  \rangle; y_i  \right)  x_{ij} \sim \subg(8n L^2 M^2), \forall j $.
\end{lemma}
In particular, if $|x_{i,j} | \le M, \ \forall i,j$, and if Assumption \ref{asu1} holds, then Hoeffding's lemma guarantees that $\partial f\left( \langle \mathbf{x}_i,  \B{\beta}^*  \rangle; y_i  \right) x_{i,j} \sim \subg(L^2M^2), \ \forall i,j$. Thus Assumption  \ref{asuSG} is satisfied. Assumption  \ref{asuSG} is also satisfied if the observations $\mathbf{x}_1,\ldots,\mathbf{x}_n$ are independently drawn from a multivariate centered Gaussian distribution. Hence, Assumption \ref{asuSG} is rather mild. It is considerably much weaker than Assumption (A1) by \cite{L1-SVM} which imposes a finite bound on the L2 norm of each column of $\mathbb{X}$.

\subsection{Restricted eigenvalue conditions}

The next assumption draws inspiration from the restricted eigenvalue conditions defined for regression problems \citep{slope, lasso-dantzig}. In particular, for an integer $k$, Assumption \ref{asu4}$.1$ ensures that some random variable is upper-bounded on the set of $k$ sparse vectors.  Similarly, Assumption \ref{asu4}$.2$ ensures that the quadratic form associated to the Hessian matrix $\nabla^2 \mathcal{L}(\B{\beta}^*) $ is lower-bounded on a cone of $\mathbb{R}^{p}$.
\begin{asu} \label{asu4}
	Let $k \in \left\{1, \ldots ,p \right\}$. Assumption \ref{asu4}$.1(k)$ is satisfied if there exists a nonnegative constant $\mu(k) $ such that almost surely:
	$$ \mu(k) \ge \sup \limits_{  \mathbf{z} \in \mathbb{R}^{p} : \ \| \mathbf{z}  \|_0 \le k  } \frac{ \sqrt{k} \| \mathbb{X} \mathbf{z}  \|_1  }{ \sqrt{n} \|\mathbf{z}\|_1  } >0.
	$$
	Let $\gamma_1, \gamma_2$ be two non-negative constants. Assumption \ref{asu4}$.2(k, \gamma)$ holds if there exists a nonnegative constant $\kappa(k, \gamma_1, \gamma_2)$ which almost surely satisfies:
	$$ 0 < \kappa(k, \gamma_1, \gamma_2) \le \inf \limits_{| S | \le k } \ \inf \limits_{ \substack{\mathbf{z} \in \Lambda(S, \gamma_1, \gamma_2) } } \frac{ \| \mathbf{z}^T \nabla^2 \mathcal{L}(\B{\beta}^*)  \mathbf{z}  \|_2  }{ \|\mathbf{z}\|_2  },
	$$
	where $\gamma = (\gamma_1, \gamma_2)$ and for every subset $S \subset \left\{1, \ldots ,p\right\}$, the cone $\Lambda(S, \gamma_1, \gamma_2) \subset \mathbb{R}^{p}$ is defined as:
	$$\Lambda(S, \gamma_1, \gamma_2) = \left\{ \mathbf{z} \in \mathbb{R}^{p}: \ \| \mathbf{z}_{S^c}  \|_1 \le \gamma_1 \| \mathbf{z}_{S}  \|_1 + \gamma_2 \| \mathbf{z}_{S}  \|_2 \right\}.$$
	We refer to Assumption \ref{asu4}$( k, \gamma$) when both Assumptions $4.1(k)$ and $4.2(k, \gamma)$ are assumed to hold.
\end{asu}
In the SVM framework, \citet{L1-SVM} define Assumption (A4): it is similar to our Assumption $4.2(k, \gamma)$ but it considers a different cone of $ \mathbb{R}^{p}$. In addition, their Assumption (A3) defines $\mu(k)$ as an upper bound of the quadratic form associated to $n^{-1/2} \mathbb{X}^T \mathbb{X}$ -- restricted to the set of $k$ sparse vectors. That is, under their definition, $\|\mathbb{X} \mathbf{z}  \|_2 / \sqrt{n} \le \mu(k) \|\mathbf{z}\|_2, \ \forall \mathbf{z} : \| \mathbf{z}  \|_0 \le k$. Our Assumption \ref{asu4}$.1(k)$ is stronger: when satisfied, we can recover Assumption (A3) since that 
\begin{align*}
\begin{split}
& \forall \mathbf{z} \in \mathbb{R}^{p} : \ \| \mathbf{z}  \|_0 \le k,\\
&\| \mathbb{X} \mathbf{z}  \|_2 / \sqrt{n} \le \| \mathbb{X} \mathbf{z}  \|_1 / \sqrt{n} \le \mu(k) \|\mathbf{z}\|_1 / \sqrt{k} \le \mu(k) \|\mathbf{z}\|_2 
\end{split}
\end{align*}
where we have used Cauchy-Schwartz inequality on the $k$ sparse vector $\mathbf{z}$. However, Assumption \ref{asu4}$.1(k)$ uses an L1 norm, more naturally associated to the class of L1-regularized estimators studied in this work.

\medskip

Similarly, in the Logistic Regression case \cite{Wainwright-logreg} consider a dependency and incoherence conditions for the population Fisher information matrix (Assumptions A1 and A2). Finally, Assumption D.4 for Quantile Regression \citep{quantile-reg} is a uniform Restricted Eigenvalue condition.

\subsection{Growth condition}

Since $\B{\beta}^*$ minimizes the theoretical loss, it holds  $\nabla L(\B{\beta}^*) = 0$. In particular, under Assumption \ref{asu4}$.2(k^*, \gamma)$, the theoretical loss evaluated on the family of cones $\Lambda(S, \gamma_1, \gamma_2)$ -- where $|S| \le k^*$ -- is lower-bounded by a quadratic form around $\B{\beta}^*$. By continuity, we define the maximal radius on which the following lower-bound holds:
\begin{align*}
\begin{split}
r(k^*) = \max \left\{  
r: \begin{array}{ll} \mathcal{L}(\B{\beta}^* + \mathbf{z} ) \ge \mathcal{L}(\B{\beta}^*) + \frac{\kappa(k^*)}{4}  \| \mathbf{z} \|_2^2\\
\forall  S \subset (p):  |S| \le k^*, \ \\
\forall \mathbf{z} \in \Lambda(S): \| \mathbf{z} \|_1 \le r
\end{array}\right\}
\end{split}
\end{align*}
where the notations $r(k^*), \ \kappa(k^*)$ and $\Lambda(S)$   are shorthands for $r(k^*, \gamma_1, \gamma_2), \ \kappa(k^*, \gamma_1, \gamma_2)$ and  $\Lambda(S, \gamma_1, \gamma_2)$. This definition is similar to the one proposed by \citet{quantile-reg} in the proof of Lemma $(3.7)$ . We now define a growth condition which gives a relation between the number of samples $n$, the dimension space $p$, our constants introduced in Assumption \ref{asu4}, and a parameter $\delta$. 
\begin{asu} \label{asu5}
	Let $\delta \in (0,1)$. We say Assumption \ref{asu5}$.1(k^*)$ is satisfied if $p \le k^* \sqrt{k^*}$. In addition, Assumption \ref{asu5}$.2( k^*, \gamma, \delta$) is said to hold if the parameters $n, \ p, \ k^*$ satisfy:  
	\begin{align*}
	\begin{split}
	\frac{\kappa(k^*)}{ 16 \alpha L  }  r(k^*) \ge  & 
	3M  \sqrt{ \frac{ k^* \log\left( 2pe/k^* \right) \log\left( 2/\delta \right) }{n} }
	+ 7 \mu(k^*)  \sqrt{ \frac{ \log(3) + \log\left(p /k^* \right) / k^* + \log\left( 2/ \delta \right) }{n}}.
	\end{split}
	\end{align*}
	We refer to Assumption \ref{asu5}$( k^*,  \gamma, \delta$) when both Assumptions $5.1(k^*)$ and $5.2(k^*, \gamma, \delta)$ hold.
\end{asu}
Assumption \ref{asu5} is similar to Equation (17) from \cite{Wainwright-logreg} for Logistic Regression. \cite{quantile-reg} also require a growth condition for Theorem 2 to hold for Quantile Regression. Consequently, as we discussed, Assumptions \ref{asu1}-\ref{asu5} are common assumptions or similar to existing ones in the literature. The next section uses our framework to derive upper bounds for L2 coefficients estimation scaling with the parameters $n, \ p, \  k^*$.


\section{Main results}\label{sec: error-bound}
This section establishes the following theorem:
\begin{theorem} \label{main-results}
	Let $\delta \in \left(0, 1 \right)$, $\alpha>1$ and assume Assumptions \ref{asu1}-\ref{asuSG}, \ref{asu4}($k^*, \gamma)$ and   \ref{asu5}($ k^*,  \gamma, \delta)$ hold -- where $\gamma=(\gamma_1, \gamma_2)$ and  $\gamma_1:= \tfrac{\alpha}{\alpha -1}$, $\gamma_2:= \tfrac{\sqrt{k^*}}{\alpha -1}$.
	
	Then, the empirical estimator $ \hat{\B{\beta}}$, defined as a solution of Problem \eqref{learning} for the  regularization parameter $\lambda = 12 \alpha L M \sqrt{ \frac{ \log(2pe/k^*) }{n} \log(2/ \delta)  }$, satisfies with probability at least $1-\delta$:
	\begin{align}\label{bound-expectation}
	\begin{split}	
	\| \hat{\B{\beta}} - \B{\beta}^*\|_2 \lesssim &
	 \frac{ \alpha L M }{\kappa(k^*)} \sqrt{ \frac{ k^* \log\left( p/k^* \right) \log\left( 2/\delta \right) }{n} } + \frac{ \alpha L \mu(k^*)  }{ \kappa(k^*)} \sqrt{ \frac{ \log(3) + \log\left(p /k^* \right) / k^* + \log\left( 2/ \delta \right) }{n}}.
	\end{split}
	\end{align}
\end{theorem}
This upper bound scales as $\left((k^*/n) \log(p/k^*)\right)^{1/2}$. It strictly improves over existing results. Note that our estimator is not adaptative to unknown sparsity: the regularization parameter $\lambda$ depends upon $k^*$. 

\medskip

The proof of Theorem \ref{main-results} is presented in Appendix \ref{sec: appendix_main-results}. It relies on two essential steps: a cone condition and a restricted strong convexity condition: these results are respectively derived in Theorems \ref{cone-condition} and \ref{restricted-strong-convexity}. The two terms of the sum in Equation \eqref{bound-expectation} are related to the two parameters $\lambda$ and $\tau$ introduced and fixed respectively in these theorems. 

\medskip

In addition, Theorem \ref{main-results} holds for any $\delta \le 1$. Thus, we obtain by integration the following bound in expectation. The proof is presented in Appendix \ref{sec: appendix_main-corollary}.
\begin{cor} \label{main-corollary}
	If the assumptions presented in Theorem are satisfied for a small enough $\delta$, then:
	\begin{equation*}
	\mathbb{E} \| \hat{\B{\beta}} - \B{\beta}^*\|_2 \lesssim \frac{\alpha L}{\kappa(k^*)} ( \mu(k^*) +M  ) \sqrt{ \frac{k^* \log\left( p /k^* \right)}{n} }.
	\end{equation*}
\end{cor}

The rest of this section follows through the steps required to prove Theorem \ref{main-results} and Corollary  \ref{main-corollary}.

\subsection{Cone condition}\label{sec: cone condition}

Similarly to the regression case~\citep{lasso-dantzig, slope}, we first derive a cone condition which applies to the difference between the empirical and theoretical minimizers. That is, by selecting a suitable regularization parameter, we show that this difference belongs to the family of cones $\Lambda(S, \gamma_1, \gamma_2)$ of $\mathbb{R}^{p}$ defined in Assumption \ref{asu4}.

\begin{theorem}\label{cone-condition}
	Let $\delta \in \left(0, 1 \right)$ and assume that Assumptions \ref{asu1} and \ref{asuSG} are satisfied. Let $\alpha\ge 2$.
	
	Let $\hat{\B{\beta}}$ be a solution of Problem  \eqref{learning} with parameter $\lambda = 12 \alpha L M \sqrt{  \frac{ \log(2pe/k^*) }{n} \log(2/ \delta)  }$.
	Then it holds with probability at least $1 - \frac{\delta}{2}$:
	
	$$\mathbf{h}:= \hat{\B{\beta}} - \B{\beta}^* \in \Lambda \left( S_0, \ \gamma_1:= \frac{\alpha}{\alpha -1},  \ \gamma_2:= \frac{\sqrt{k^*}}{\alpha -1} \right),$$
	where $S_0$ is the subset of indices of the $k^*$ highest coefficients of $\mathbf{h}$.
\end{theorem}

The regularization parameter $\lambda$ is selected so that it dominates the sub-gradient of the loss $f$ evaluated at the theoretical minimizer $\B{\beta}^*$. The proof is presented in Appendix \ref{sec: appendix_cone-condition}: it uses a new result to control the maximum of independent sub-Gaussian random variables. As a result, our cone condition is stronger than the ones proposed by \cite{L1-SVM} and \cite{Wainwright-logreg}: their value of $\lambda^2$ is of the order of $(k^*/n)\log(p)$ whereas ours scales as $(k^*/n) \log(p/k^*)$.

\subsection{A supremum result }\label{sec: supremum}

The next Theorem \ref{hoeffding-sup} is an essential step to obtain our main Theorem \ref{main-results}. It derives a control of the supremum of the difference between an empirical random variable and its expectation. This supremum is controlled over a bounded set of sequences of $k$ sparse vectors with disjoint supports. The restricted strong convexity condition derived in Theorem \ref{restricted-strong-convexity} is a consequence of Theorem \ref{hoeffding-sup}. 

\medskip

To motivate this theorem, it helps considering the difference between the usual regression framework and our framework for classification problems. The linear regression case assumes the generative model $\B{y} = \B{X} \B{\beta}^* + \B{\epsilon}$. Therefore, with the notations of Theorem \ref{hoeffding-sup},  $\Delta(\B{\beta}^*, \B{z}) = \frac{1}{n} \| \B{X}\B{z} \|_2^2 - \frac{2}{n} \B{\epsilon}^T\B{X}\B{z} $. By combining a cone condition (similar to Theorem 1) with an upper-bound of the term $\epsilon^T \B{X} \B{z} $, we can obtain a restricted strong convexity similar to Theorem \ref{restricted-strong-convexity}. However, in the classification case, $\B{\beta}^*$ is defined as the minimizer of the theoretical risk. Two majors differences appear: (i) we cannot simplify $\Delta( \B{\beta}^*, \B{z} )$ with basic algebra, (ii) we need to introduce the expectation $\mathbb{E} (\Delta(\B{\beta}^*, \B{z} ))$ and to control the quantity $ | \mathbb{E}(\Delta(\B{\beta}^*,  \B{z} )) - \Delta(\B{\beta}^*, \B{z}) |$. Theorem \ref{hoeffding-sup} helps expliciting the cost to pay for this control.

\begin{theorem} \label{hoeffding-sup} We define $\forall  \mathbf{w}, \mathbf{z} \in \mathbb{R}^p$:
	$$\Delta(\mathbf{w},  \mathbf{z}) = \frac{1}{n} \sum_{i=1}^n f \left( \langle \mathbf{x}_i,  \mathbf{w}+ \mathbf{z}  \rangle ;  y_i \right)  - \frac{1}{n} \sum_{i=1}^n  f \left( \langle \mathbf{x}_i,  \mathbf{w} \rangle ;  y_i \right).$$ 
	Let $k\in \left\{1, \ldots, p\right\}$ and  $S_1, \ldots S_q$ be a partition of $\left\{1, \ldots, p\right\}$ with $q = \lceil p/k \rceil$ and $|S_j| \le k, \forall j$. 
	
	Let $\tau(k) = 14 L  \mu(k) \sqrt{\frac{\log(3) }{n} + \frac{\log\left( 4 p/k \right) }{nk}   +  \frac{ \log\left( 2/\delta\right)  }{nk}  } $ and assume that Assumptions \ref{asu1}, \ref{asu4}$.1(k)$ and \ref{asu5}$.1(k$) hold. Then, for any $\delta \in (0,1)$, it holds with probability at least $1-\frac{\delta}{2}$:
	$$  \sup \limits_{ \substack{  \mathbf{z}_{S_1}, \ldots, \mathbf{z}_{S_q}  \in \mathbb{R}^{p}: \\  \Supp(\mathbf{z}_{S_j})  \subset S_j \ \forall j \\ \ \ \| \mathbf{z}_{S_j}  \|_1 \le 3R \ \ \ \forall j } }   \left\{ 
	\sup \limits_{\ell=1,\ldots,q} 
	\left\{ \Omega \left( \mathbf{w}_{\ell-1}, \ \mathbf{z}_{S_\ell} \right)
	\right\} \right\}
	\le 0, \textnormal{with}$$
	\begin{align*}
	\begin{split}
	\Omega \left( \mathbf{w}_{\ell-1},\mathbf{z}_{S_\ell} \right) 
	:=& \left|  \Delta \left( \mathbf{w}_{\ell-1}, \mathbf{z}_{S_\ell} \right)   - \mathbb{E}\left(  \Delta \left( \mathbf{w}_{\ell-1}, \mathbf{z}_{S_\ell} \right) \right) \right| 
	-  \tau(k) \| \mathbf{z}_{S_\ell}\|_1.
	\end{split}
	\end{align*} 
	$\Supp(.)$ refers to the support of a vector and we define $\mathbf{w}_{\ell} = \B{\beta}^*  + \sum \limits_{j=1}^{\ell} \mathbf{z}_{S_j}, \forall \ell$. 
\end{theorem}
The proof is presented in Appendix \ref{sec: hoeffding-sup}. It uses Hoeffding's inequality to obtain an upper bound of the inner supremum for any sequence of $k$ sparse vectors. The result is extended to the outer supremum with an $\epsilon$-net argument.

\subsection{Restricted strong convexity condition}\label{sec: restricted}

Theorem \ref{hoeffding-sup} applies to a sequence of $k$ sparse vectors with disjoint supports. In particular we can fix $k=k^*$ and consider $\mathbf{h}=\hat{\B{\beta}} - \B{\beta}^*$. In addition, we can exploit the minimality of $\B{\beta}^*$ and the cone condition proved in Theorem \ref{cone-condition}. By pairing these points, we derive the next Theorem \ref{restricted-strong-convexity}. It says that the loss $f$ satisfies a restricted strong convexity \citep{M-estimators} with curvature $ \kappa(k^*) /4$ and L1 tolerance function.

\begin{theorem} \label{restricted-strong-convexity}Let $\textbf{h} = \hat{\B{\beta}} - \B{\beta}^*$ and $\delta \in (0,1)$. Under the notations of Theorem \ref{hoeffding-sup}, if Assumptions \ref{asu1}-\ref{asuSG}, \ref{asu4}$(k^*, \gamma)$ and \ref{asu5}($ k^*,  \gamma, \delta)$ are satisfied, then it holds with probability  at least $1 - \delta$:
	$$\Delta\left(\B{\beta}^*, \mathbf{h} \right) \ge \frac{1}{4}  \kappa(k^*) \left\{  \|\mathbf{h}\|_2^2 \wedge  r(k^*) \|\mathbf{h}\|_2  \right\}  - \tau(k^*) \| \mathbf{h} \|_1.$$
\end{theorem}
The proof is presented in Appendix \ref{sec: appendix_restricted-strong-convexity}. Let us note that our parameter  $\tau(k^*)^2$ scales  as $n^{-1}(1 + \log(p/k^*)/k^*)$, whereas \cite{L1-SVM}, \cite{Wainwright-logreg} and \cite{M-estimators} all propose a parameter scaling as $n^{-1}k^* \log(p)$. Our restricted strong convexity condition is stronger. We later use the cone condition derived in Theorem 1 to convert the L1 tolerance function into the L2 norm used for coefficients estimation.

\subsection{Deriving Theorem \ref{main-results} and Corollary \ref{main-corollary}  }

Our main bounds -- presented in Theorem \ref{main-results} and Corollary \ref{main-corollary} -- follow from  the two preceding Theorems  \ref{cone-condition} and \ref{restricted-strong-convexity}. The proofs are respectively presented in Appendix \ref{sec: appendix_main-results} and \ref{sec: appendix_main-corollary}. Our family of L1-regularized L1-constrained estimators reach a bound that strictly improve over existing results. Our rate is the best known for the classification problems considered here, and it holds both with high probability and in expectation.

\section{Algorithm and upper bounds for Slope estimator}\label{sec: slope}

This section introduces the Slope estimator -- originally presented for the linear regression case \citep{slope-introduction, slope-proximal} --  to our class of problems. We propose a tractable algorithm to compute it and study its statistical performance.

\subsection{Introducing Slope for classification}

We consider a sequence $\lambda \in \mathbb{R}^p$ such that $\lambda_1 \ge \ldots \ge \lambda_p >0$, and we note $\mathcal{S}_p$ the set of permutations of $\left\{1,\ldots ,p \right\}$. The Slope regularization is defined as:
\begin{equation} \label{slope-norm}
| \B{\beta} |_S = \max \limits_{\phi \in \mathcal{S}_p } \sum_{j=1}^p | \lambda_j |  | \beta_{\phi(j)} | = \sum_{j=1}^p \lambda_j | \beta_{(j)} |,
\end{equation}
where $|\beta_{(1)}| \ge \ldots \ge |\beta_{(p)}| $ is a non-increasing rearrangement of $\B{\beta}$. Consequently for $\eta>0$, we define the Slope estimator $\B{\hat{\beta}}$ as the solution of the convex minimization problem:
\begin{equation} \label{slope}
\min \limits_{ \B{\beta} \in \mathbb{R}^{p}  }  \frac{1}{n} \sum_{i=1}^n  f \left( \langle \mathbf{x_i},  \B{\beta} \rangle ;  y_i \right)  + \eta | \B{\beta} |_S.
\end{equation}
The approach presented herein uses a proximal gradient algorithm -- with Nesterov smoothing \citep{smoothing-math} in the case of the hinge loss and quantile regression loss -- to solve Problem \eqref{slope}, extending the original definition of Slope \citep{slope-introduction} to a larger class of loss functions.  Recently, \citet{dedieu2019} combined similar first order methods for non-smooth convex optimization with column with constraint generation algorithms to solve linear SVM with sparsity-inducing regularization when the number of samples and / or features is of the order of hundreds of thousands.

\subsection{Smoothing the hinge loss  } \label{sec: smoothing-hinge}

The method described in Section \ref{sec: thresholding} to solve Problem \eqref{slope} requires $f(.,y)$ to be differentiable with Lipschitz-continuous gradient. Among the loss functions considered in Section \ref{sec: framework}, only the  logistic regression loss satisfies this condition.

\medskip

To handle the non-smooth hinge loss, we use the smoothing scheme pioneered by \cite{smoothing-math}. We construct a convex function $g^{\tau}$ with continuous Lipschitz gradient, which approximates the hinge loss for $\tau \approx 0$. 
Let us first note that $\max(0,x) = \frac{1}{2}(x + |x|) =  \max_{|w| \le 1} \frac{1}{2}(x + wx)$ as this maximum is achieved for $\sign(x)$. Consequently the hinge loss can be expressed as a maximum over the L$_\infty$ unit ball:
$$\frac{1}{n}\sum \limits_{i=1}^n \max(z_i, 0) =  \max \limits_{\|w\|_{\infty} \le 1}  \frac{1}{2n}  \sum \limits_{i=1}^n \left[  z_i  + w_i z_i  \right], $$
where $z_i = 1 - y_i \mathbf{x}_i^T\B{\beta}, \ \forall i$. We apply the technique suggested by  \cite{smoothing-math} and define for $ \tau >0$ the smoothed hinge loss:
\begin{align}\label{smooth-hinge-def}
&g^{\tau}(\B{\beta}) =  \max \limits_{\|w\|_{\infty} \le 1} \frac{1}{2n}  \sum \limits_{i=1}^n  \left[  z_i   + w_i z_i  \right] - \frac{\tau}{2n} \|w\|_2^2.
\end{align}
Let $\textbf{w}^{\tau}(\B{\beta}) \in \mathbb{R}^n:  \ w^{\tau}_i(\B{\beta}) = \min\left( 1, \frac{1}{2\tau} | z_i | \right) \sign(z_i), \ \forall i$ be the optimal solution of the right-hand side of Equation \eqref{smooth-hinge-def}. The gradient of $g^{\tau}$ is expressed as:
\begin{equation}
\nabla g^{\tau}( \B{\beta} ) = - \frac{1}{2n} \sum \limits_{i=1}^{n}(1+w_i^{\tau}(\B{\beta}) )y_i \mathbf{x}_i \in \mathbb{R}^{p},
\end{equation}
and its associated Lipschitz constant is derived from the next theorem.
\begin{theorem} \label{lipschitz}
	Let $\mu_{\max} (n^{-1}\mathbb{X}^T \mathbb{X})$ be the highest eigenvalue of $n^{-1}\mathbb{X}^T \mathbb{X}$. Then $\nabla g^{\tau}$ is Lipschitz continuous with constant $ C^{\tau} = \mu_{\max} (n^{-1}\mathbb{X}^T \mathbb{X}) / 4 \tau$.
\end{theorem}
The proof is presented in Appendix  \ref{sec: appendix_lipschitz}. It follows \cite{smoothing-math} and uses first order necessary conditions for optimality. We mention how to adapt the theorem to the quantile regression loss.

\subsection{Thresholding operator for Slope}\label{sec: thresholding}
We note $g(\B{\beta}) = \frac{1}{n} \sum_{i=1}^n  f \left( \langle \mathbf{x_i},  \B{\beta} \rangle ;  y_i \right)$. Problem \eqref{slope} can be equivalently formulated as:
\begin{equation}\label{smooth-slope} 
\min \limits_{ \B{\beta} \in \mathbb{R}^{p}  } g(\B{\beta})    + \eta | \B{\beta} |_S, 
\end{equation}
We now require $g$ to be a differentiable loss with $C$-Lipschitz continuous gradient. When $f$ is the hinge or quantile regression loss we replace $g$ with $g^{\tau}$ as defined in Section \ref{sec: smoothing-hinge}.  For $D\ge C$, we upper-bound $g$ around any $\B{\alpha} \in \mathbb{R}^{p}$ with the quadratic form $Q_D(\B{\alpha},.)$ defined as the right-hand side of the equation:
\begin{equation}\label{convex-continuous-gradient}
g( \B{\beta} ) \le g(\B{\alpha}) + \nabla g(\B{\alpha})^T(\B{\beta}-\B{\alpha}) + \frac{D}{2} \| \B{\beta} - \B{\alpha} \|_2^2.
\end{equation}
We approximate the solution of Problem $\eqref{slope}$  by considering the loss $Q_D$ and solving the problem:
\begin{equation} \label{key-proxy-slope}
\begin{split}	
\argmin_{ \B{\beta}  } Q_D(\B{\alpha},\B{\beta}) + \eta | \B{\beta} |_S 
&= \argmin_{ \B{\beta}  } \frac{1}{2}  \left\| \B{\beta} - \left(\B{\alpha} - \frac{1}{D} \nabla g(\B{\alpha})  \right) \right\|_2^2 + \frac{\eta}{D}  | \B{\beta} |_S\\
&= \argmin_{ \B{\beta}  } \frac{1}{2}  \left\| \B{\beta} - \B{\gamma} \right\|_2^2 +  \sum_{j=1}^p \tilde{\eta}_j | \beta_{(j)} |,
\end{split}	
\end{equation}
where $\B{\gamma} = \B{\alpha} - \frac{1}{D} \nabla g(\B{\alpha})$ and $\tilde{\eta}_j = \frac{\eta}{D} \lambda_j, \ \forall j$. To solve Problem \eqref{key-proxy-slope}, we need to derive the proximal operator of the sorted L1 norm. The next Lemma \ref{lemma-thresholding-slope} does so by noting that the signs of the quantities $\beta_j$ and $\gamma_j$ are all identical.
\begin{lemma}  \label{lemma-thresholding-slope}
	Let us assume that $\tilde\gamma_1\ge \ldots \ge \tilde\gamma_p \geq 0$. Since $\tilde{\eta}_1\ge \ldots \ge \tilde{\eta}_p \geq 0$, the solution of Problem \eqref{key-proxy-slope} can be derived from the solution of the problem:
	\begin{equation}\label{slope-thresholded}
	\begin{myarray}[1.1]{c c c r}
	\min \limits_{\B{\beta} \in \mathbb{R}^{p } } & \frac{1}{2}  \left\| \B{\beta} - \B{\tilde{\gamma}} \right\|_2^2 +  \sum \limits_{j=1}^p \tilde{\eta}_j  \beta_j & \\
	\sbt 
	& \;\;\;\;  \beta_1 \ge \ldots \ge \beta_p \ge 0. & 
	\end{myarray}
	\end{equation}
\end{lemma} 
\cite{slope-proximal} proposed an efficient proximal algorithm to solve Problem \eqref{slope-thresholded} called FastProxSL1: it is guaranteed to terminate in at most $p$ iterations. We denote by $\mathcal{T}_{ \left\{ \tilde{\eta}_j \right\}  }( \B{\gamma})$ a solution to Problem \eqref{key-proxy-slope}.

\subsection{First order algorithm}\label{sec :AGD} 
The following algorithm applies the accelerated gradient descent method \citep{FISTA}  on the smoothed version of the Slope Problem \eqref{smooth-slope} by using the above thresholding operator. 
The iterations continue till the algorithm converges or a  maximum number of iterations $T_{\max}$ is reached.
\medskip

\textbf{Input: }$\mathbf{X}$, $\mathbf{y}$, a sequence of Slope coefficients  $\left\{ \lambda_j \right\} $, a regularization parameter $\eta$,  a stopping criterion $\epsilon$, a maximum number of iterations $T_{\max}$.
\textbf{Output: } An approximate solution $\B{\beta}$ for the smoothed Slope Problem \eqref{smooth-slope}.
\begin{enumerate}
	\item Initialize $T=1$, $q_1 = 1$,  $\B{\beta}_1 = \B{\delta}_0=0$.
	\item : While $\| \B{\beta}_{T} - \B{\beta}_{T-1} \|_2^2  > \epsilon $ and $T<T_{\max}$ do:
	\begin{enumerate}
		\item Compute $\B{\delta}_{T} = \mathcal{T}_{ \left\{ \eta \lambda_j / C  \right\}  } \left( \B{\beta}_T - \frac{1}{C } \nabla g (\B{\beta }_T )  \right)  $.
		\item Define $q_{T+1} = \frac{1 + \sqrt{1 + 4 q_T ^2} }{2}$ and compute  $\B{\beta}_{T+1} = \B{\delta}_{T} + \frac{q_T - 1}{q_{T+1}} (\B{\delta}_{T} - \B{\delta}_{T-1})$.
	\end{enumerate}
\end{enumerate}

\subsection{Error bounds for Slope}

We extend our previous case and study under our framework the theoretical properties of a Slope estimator. In particular, we consider the L1-constrained Slope estimator:
\begin{equation} \label{slope-constrained}
\min \limits_{ \B{\beta} \in \mathbb{R}^{p}: \ \| \B{\beta} \|_1 \le 2R }  \frac{1}{n} \sum_{i=1}^n  f \left( \langle \mathbf{x_i},  \B{\beta} \rangle ;  y_i \right)  + \eta | \B{\beta} |_S.
\end{equation}
The study of Slope share a lot of similarities with our previous work for L1-regularized estimators. First, we derive the following cone condition:
\begin{theorem} \label{cone-condition-slope}
	Let $\delta \in \left(0, 1 \right)$ and $\alpha \ge 2$. We fix the Slope coefficients $\lambda_j =  \sqrt{  \log(2pe/j) }, \forall j$, and assume Assumptions \ref{asu1} and \ref{asuSG} hold. Then the Slope estimator defined as a solution of Problem \eqref{slope-constrained} for the regularization parameter $\eta = 14 \alpha L M \sqrt{ n^{-1} \log(6/ \delta)}$ satisfies with probability at least $1-\frac{\delta}{2}:$
	$$\hat{\B{\beta}}  - \B{\beta}^*\in \Gamma \left( k^*, \omega^*= \frac{\alpha +1}{\alpha -1} \right),$$
	where for every $k \in \left\{1, \ldots ,p \right\}$ and $\omega >0$, the cone $\Gamma(k, \omega) $ is defined as:
	$$\Gamma(k, \omega) = \left\{ \mathbf{z} \in \mathbb{R}^{p}: \ \sum_{j=k+1}^p \lambda_j |z_{(j)}| \le   \omega \sum_{j=1}^{k} \lambda_j |z_{(j)}| \right\}$$
	with $| z_{(1)}| \ge \ldots \ge |z_{(p)}|, \ \forall \B{z}$.
\end{theorem}
The proof is presented in Appendix \ref{sec: cone-condition-slope}. We consequently adapt Assumption \ref{asu4bis}$.2$ to the new family of cones $\Gamma(k, \omega)$ introduced in Theorem \ref{cone-condition-slope}.
\begin{asu}\label{asu4bis}
	Let  $k \in \left\{1, \ldots ,p \right\}$ and $\omega >0$. Assumption \ref{asu4bis}$.2(k, \omega)$ is said to hold if there exists a nonnegative constant $\kappa(k, \omega)$ such that:
	$$ 0 < \kappa(k, \omega) \le \ \inf \limits_{ \mathbf{z} \in \Gamma(k,\omega)  }  \frac{ \| \mathbf{z}^T \nabla^2\mathcal{L}(\B{\beta}^*)  \mathbf{z}  \|_2  }{ \|\mathbf{z}\|_2  }.
	$$
\end{asu}
Similarly, we define a new growth condition -- Assumption $8(k, \omega, \delta)$ -- which adapts Assumption~6 to Slope by replacing $\kappa(k)$ with $\kappa(k, \omega)$ defined above. The following result holds for Slope. The proof is presented in Appendix \ref{sec: results-slope}.
\begin{cor} \label{main-results-slope}
	Assume Assumptions \ref{asu1}-\ref{asuSG}, Asumptions \ref{asu4bis}$(k^*, \omega^*)$ and   $8(k^*, \omega^*, \delta)$ hold for a small enough $\delta$, where $\omega^*$ is defined in Theorem \ref{cone-condition-slope}.
	
	Then the bounds presented in Theorem \ref{main-results} and Corollary \ref{main-results} are achieved by a Slope  estimator, defined as a solution of Problem \eqref{slope} for the coefficients $\lambda_j =  \sqrt{  \log(2pe/j) }, \forall j$ and the regularization parameter $\eta = 14 \alpha L M n^{-1}\sqrt{ \log(6/ \delta)}$ -- where $\alpha\ge2$.
\end{cor}

This Slope estimator adapts to unknown sparsity while achieving the same bound than the L1-regularized estimator studied in Theorem \ref{main-results} and Corollary \ref{main-results}.

\subsection{Simulations}

We finally compute a family of Slope estimators and demonstrate its empirical performance -- for L2 coefficients estimations and misclassification accuracy -- we compare to L1 and L2-regularized estimators.

\medskip 

\noindent {\textbf{Data Generation:}} We consider $n$ independent realizations of a $p$ dimensional multivariate normal centered distribution, with only $k^*$ dimensions being relevant for classification. Half of the samples are from the $+1$ class and have mean $\mu_+ = (\mathbf{1}_{k^*}, \ \mathbf{0}_{p-k^*} )$. The other half are from the $-1$ class and have mean $\mu_- = - \mu_+$. We consider a covariance matrix $\Sigma_{ij} = \rho$ if $i \ne j$ and $1$ otherwise. The data of both $\pm1$ classes respectively have the distribution: $\forall i, \  x_i^{\pm} \sim \mathbb{N}(\mu_{\pm},\Sigma)$. 

\medskip

\noindent {\textbf{Competitors:}} Table \ref{table:results} compares the performance of 3 approaches -- each associated to a different regularization -- for both the SVM and the Logistic Regression problems. Method \textbf{(a)} computes a family of L1-regularized estimators for a decreasing geometric sequence of regularization parameters $\eta_0> \ldots > \eta_M$. We start from a high enough $\eta_{0}$  so that the solution of Problem \eqref{learning} is the $0$ estimator and we fix $\eta_M < 10^{-4} \eta_0$. For the hinge loss, we solve the Linear Programming L1-SVM problem with the commercial LP solver \textsc{Gurobi} version $6.5$ with Python interface. The L1-regularized Logistic Regression is solved with \textsc{scikit-learn} Python package.  In addition, method \textbf{(b)} returns a family of L2-regularized estimators with \textsc{scikit-learn} package: we start from $\eta_0 = \max_i \left\{ \| \mathbf{x_i} \|_2^2 \right\}$ as suggested by \cite{warm-start}.  Finally, method \textbf{(c)} computes a family of Slope-regularized estimators, using the first order algorithm presented in Section \ref{sec :AGD} for $\tau= 0.2$. The Slope coefficients $\left\{ \lambda_j \right\}$ are the ones proposed in Theorem \ref{cone-condition-slope}; the set of parameters $\left\{\eta_i\right\}$ is identical to method \textbf{(a)}.

\medskip

\noindent {\textbf{Metrics:} Following our theoretical results, we want to find the estimator which minimizes the L2 estimation error:  
	$$\left\| \frac{\hat{\B{\beta}}  }{\| \hat{\B{\beta}} \|_2 } - \frac{\B{\beta}^*  }{\| \B{\beta}^*  \|_2 }   \right\|_2,$$
	where $\B{\beta}^* $ is the theoretical minimizer. $\B{\beta}^* $ is computed on a large test set with $10,000$ samples: we solve the SVM or Logistic Regression problem with a very small regularization coefficient on the $k^*$ columns relevant for classification. We also study the misclassification performances on this same test set. For each family returned by the methods \textbf{(a)}, \textbf{(b)} and \textbf{(c)}, we only keep the estimator with lowest misclassification error on an independent validation set of size $10,000$.  
	
	Table \ref{table:results} compares the L2 estimation error (\textbf{L2-E}), and the test misclassification error (\textbf{Misc}) -- of these 3 estimators  selected on the validation set. The results are averaged over 10 simulations. 

\begin{table*}[!h]
	\centering
	\caption{{\small{Averaged L2 estimation (\textbf{L2-E}) and test misclassification error (\textbf{Misc}) for the methods \textbf{(a)}, \textbf{(b)} and \textbf{(c)} over $10$ repetitions. We use varying $n,p$ values with $k^*=n/10$ and $\rho=0.1$. The Slope estimator shows impressive gains for estimating the theoretical minimizer $\B{\beta}^* $, while achieving the lowest misclassification errors.}}}
	\medskip
	\label{table:results}
	\vspace{-.5em}
	\begin{tabular}{lllllllll}
		\cmidrule{2-9} 
		&\multicolumn{2}{c}{$n=100, p=1k$} &\multicolumn{2}{c}{$n=100, p=10k$} &\multicolumn{2}{c}{$n=1k, p=1k$} & \multicolumn{2}{c}{$n=1k, p=10k$}  \\
		\cmidrule(r){2-9} 
		& L2-E  & Misc(\%)  & L2-E & Misc(\%)  & L2-E  & Misc(\%)  & L2-E  & Misc(\%)  \\
		\toprule
		L1 SVM  &  0.57 & 1.67 & 0.52 & 1.54 & 1.12  & 1.17 & 1.01 & 0.15  \\			 
		L2 SVM  &  0.54 &  1.73 & 0.52 & 1.54 & 1.11 & 0.18 & 0.91 & 0.11   \\			
		Slope SVM   & 0.34 & 1.24 & 0.37 & 1.15 & 0.94 & 0.13 & 0.83 & 0.10 \\			
		\midrule
		L1 LR & 0.48 & 1.40  & 0.46 & 1.37 & 1.04  & 0.18 & 1.04 & 0.16  \\			
		L2 LR  &  0.92 & 3.2  & 1.25 & 0.18 & 0.82  & 0.12  & 0.89 & 0.16 \\			
		Slope LR & 0.22  & 1.14  & 0.18 & 1.12 & 0.81  & 0.12 & 0.82 & 0.13  \\			
		\bottomrule
	\end{tabular}
\end{table*}

\newpage
\bibliographystyle{unsrtnat}
\bibliography{arxiv}

\newpage
\begin{appendices}
\section{Usefull properties of sub-Gaussian random variables}

This section presents useful preliminary results satisfied by sub-Gaussian random variables. In particular, Lemma \ref{upper-bound-sup} provides a probabilistic upper-bound on the maximum of independent sub-Gaussian random variables.

\subsection{Preliminary results}

Under Assumption \ref{asuSG}, the random variables $\sum \limits_{i=1}^n  \partial f\left( \langle \mathbf{x}_i,  \B{\beta}^*  \rangle, y_i  \right)  x_{ij},  \ \forall j$ are sub-Gaussian.  They all consequently satisfy the next Lemma  \ref{lemma-lecture-notes}:

\begin{lemma} \label{lemma-lecture-notes}
	Let $Z \sim \subg(\sigma^2)$ for a fixed $\sigma>0$. Then for any $t>0$ it holds 
	$$\mathbb{E}\left( \exp(tZ) \right) \le e^{4 \sigma^2 t^2}. $$
	
	In addition, for any positive integer $\ell \ge 1$ we have:
	$$\mathbb{E}\left( | Z |^{\ell} \right) \le  (2\sigma^2)^{\ell/2} \ell \Gamma(\ell/2)$$
	where $\Gamma$ is the Gamma function defined as  $\Gamma(t) =  \int_{0}^{\infty} x^{t-1} e^{-x} dx, \ \forall t>0.$ 
	
	\medskip
	
	Finally, let $Y=Z^2 - \mathbb{E}(Z^2)$ then we have
	\begin{equation}\label{eq-lemma-preliminary}
	\mathbb{E} \left( \exp\left( \frac{1}{16 \sigma^2} Y \right) \right) \le \frac{3}{2},
	\end{equation}
	and as a consequence $\mathbb{E} \left( \exp\left( \frac{1}{16 \sigma^2} Z^2\right)  \right) \le 2.$
\end{lemma}

\begin{Proof}
	The two first results correspond to Lemmas 1.4 and 1.5 from \cite{lecture-notes}. 
	
	In particular $\mathbb{E}\left( | Z |^2 \right) \le  4\sigma^2$. 
	
	In addition, using the proof of Lemma 1.12 we have:
	$$\mathbb{E} \left( \exp(tY ) \right) \le 1 + 128 t^2 \sigma^4, \ \forall |t| \le \frac{1}{16\sigma^2}.$$
	Equation \eqref{eq-lemma-preliminary} holds in the particular case where $t=1/16\sigma^2.$ 
	
	The last part of the lemma combines our precedent results with the observation that $\frac{3}{2} e^{1/4} \le 2$.
\end{Proof}

\subsection{Proof of Lemma \ref{lemma-asu3}}\label{sec: appendix_lemma}
As a first consequence of Lemma  \ref{lemma-lecture-notes}, we easily derive the proof of Lemma \ref{lemma-asu3} -- stated in Section \ref{sec: asu-SG}.

\begin{Proof}
	We note $S_i = \partial f\left( \langle \mathbf{x}_i,  \B{\beta}^*  \rangle, y_i  \right), \ \forall i $. 
	
	Since $\B{\beta}^*$ minimizes the theoretical loss, we have $\mathbb{E}(S_i x_{ij})=0, \ \forall i,j$. 
	
	\smallskip
	
	We fix $M>0$ such that: $\forall t>0,$
	$$\mathbb{P}\left(  | S_i x_{i,j} | > t  \right) \le 2 \exp \left( - \frac{t^2}{2 L^2 M^2} \right) , \ \forall i,j.$$
	Then from Lemma \ref{lemma-lecture-notes} it holds:
	$$\mathbb{E}\left( \exp(t S_i x_{ij} ) \right) \le e^{4 L^2 M^2 t^2}, \ \forall t>0,  \forall i,j.$$
	As a consequence, using Lemma \ref{lemma-lecture-notes} for the independent random variables $\left(S_1 x_{1,j},\ldots, S_n x_{n,j} \right)$, it holds $\forall t>0,$
	\begin{align*}
	\begin{split}
	\mathbb{E}\left( \exp \left( \frac{t}{\sqrt{n}} \sum_{i=1}^n S_i x_{i,j} \right) \right) 
	= \prod_{i=1}^n  \mathbb{E}\left( \exp \left( \frac{t}{\sqrt{n}}  S_i  x_{ij} \right) \right) 
	\le \prod_{i=1}^n e^{4 L^2 M^2 t^2 /n}  = e^{4 L^2 M^2 t^2}.
	\end{split}
	\end{align*}
	Let $M_1 = 2 \sqrt{2} M$, then with a Chernoff bound: 
	\begin{align*}
	\begin{split}
	\mathbb{P}\left(  \frac{1}{\sqrt{n}} \sum_{i=1}^n S_i x_{i,j}  > t  \right)  
	\le \min_{s>0} \ \exp\left( \frac{M_1^2 L^2 s^2}{2} - s t \right) 
	= \exp \left( - \frac{t^2}{2 L^2 M_1^2} \right), \ \forall t>0,
	\end{split}
	\end{align*}
	which concludes the proof.
\end{Proof}

\subsection{A bound for the maximum of independent sub-Gaussian variables}

The next two technical lemmas derive a probabilistic upper-bound for the maximum of sub-Gaussian random variables. Lemma \ref{lemma-Slope} is an extension for sub-Gaussian random variables of Proposition E.1 \citep{slope}.
\begin{lemma}\label{lemma-Slope}
	Let $g_1,\ldots g_p$ be independent sub-Gaussian random variables with variance $\sigma^2$. Denote by $(g_{(1)}, \ldots, g_{(p)})$ a non-increasing rearrangement of $(|g_1|, \ldots, |g_p|)$.  Then $\forall t>0$ and $\forall j \in \left\{1, \ldots,p \right\}$:
	$$ \mathbb{P}\left( \frac{1}{j \sigma^2} \sum_{k=1}^j g_{(k)} ^2 > t \log\left( \frac{2p}{j} \right) \right) \le \left( \frac{2p}{j} \right)^{1- \frac{t}{16}}. $$
\end{lemma}

\begin{Proof}
	We first apply a Chernoff bound: 
	\begin{align*}
	\begin{split}
	\mathbb{P}\left( \frac{1}{j\sigma^2} \sum_{k=1}^j g_{(k)} ^2 > t \log\left( \frac{2p}{j} \right) \right) 
	\le \mathbb{E} \left( \exp \left( \frac{1}{16 j \sigma^2} \sum_{k=1}^j g_{(k)} ^2 \right)  \right)   \left( \frac{2p}{j} \right)^{-\frac{t}{16}}.
	\end{split}
	\end{align*}
	
	Then we use Jensen inequality to obtain
	\begin{align*}
	\begin{split}
	\mathbb{E} \left( \exp\left( \frac{1}{16 j \sigma^2} \sum_{k=1}^j g_{(k)} ^2 \right) \right)
	\le \frac{1}{ j} \sum_{k=1}^j \mathbb{E} \left( \exp\left( \frac{1}{16 \sigma^2} g_{(k)}^2 \right) \right) 
	& \le \frac{1}{ j} \sum_{k=1}^p \mathbb{E} \left( \exp\left( \frac{1}{16 \sigma^2}  g_{k}^2 \right) \right) \\
	& \le \frac{2p}{j} \text{  with Lemma \ref{lemma-lecture-notes}}.
	\end{split}
	\end{align*}
\end{Proof}

Using Lemma \ref{lemma-Slope}, we can derive the following bound holding with high probability:

\begin{lemma}\label{upper-bound-sup}
	We consider the assumptions and notations of Lemma \ref{lemma-Slope}. In addition, we define the coefficients $\lambda_j = \sqrt{ \log(2p/j) }, \ j=1,\ldots\,p$. Then for $\delta \in \left(0, \frac{1}{2} \right)$, it holds with probability at least $1-\delta$:
	$$ \sup \limits_{j=1,\ldots,p} \left\{ \frac{ g_{(j)} }{\sigma \lambda_j} \right\} \le 12 \sqrt{ \log(1 / \delta)}.  $$ 
\end{lemma}

\begin{Proof}
	We fix $\delta \in \left(0, \frac{1}{2} \right)$ and $j \in \left\{1, \ldots,p \right\}$. We upper-bound $g_{(j)}^2$ by the average of all larger variables:
	$$g_{(j)}^2 \le \frac{1}{j} \sum_{k=1}^j g_{(k)}^2. $$
	Applying Lemma \ref{lemma-Slope} gives, for $t>0$:
	\begin{align*}
	\begin{split}
	\mathbb{P}\left( \frac{ g_{(j)}^2 }{\sigma^2 \lambda_j^2}  > t  \right) &\le \mathbb{P}\left( \frac{1}{j \sigma^2} \sum_{k=1}^j g_{(k)} ^2 > t  \lambda_j^2 \right) \le \left( \frac{j}{2p} \right)^{\frac{t}{16} -1 }.
	\end{split}
	\end{align*}
	
	We fix $t = 144 \log(1/\delta)$ and use an union bound to get:
	\begin{align*}
	\begin{split}
	\mathbb{P}\left( \sup \limits_{j=1,\ldots,p} \frac{  g_{(j)}  }{\sigma \lambda_j}  > 12 \sqrt{\log(1/\delta) }  \right) 
	\le \left( \frac{1}{2p} \right)^{9 \log(1/\delta) -1 }  \sum_{j=1}^p j^{9 \log(1/\delta)  -1 }. 
	\end{split}
	\end{align*}
	Since $\delta < \frac{1}{2}$ it holds that $9 \log(1/\delta)  -1\ge 9 \log(2) -1 > 0$, then the map $t>0 \mapsto t^{9 \log(1/\delta)  -1 }$ is increasing. An integral comparison gives:
	$$\sum_{j=1}^p j^{9 \log(1/\delta)  -1 } \le \frac{1}{2} \left( p+1 \right)^{9 \log(1/\delta) } = \frac{1}{2} \delta^{-9 \log(p+1)}.$$
	
	In addition  $9 \log(1/\delta) -1 \ge 7 \log(1/\delta)$ and 
	$$\left( \frac{1}{2p} \right)^{9 \log(1/\delta) -1 } \le \left( \frac{1}{2p} \right)^{ -7 \log(\delta)  } = \delta^{ 7\log(2p)}.$$
	
	Finally, by assuming $p\ge2$, then we have $ 7\log(2p) - 9\log(p+1) > 1$ and we conclude: 
	$$ \mathbb{P}\left( \sup \limits_{j=1,\ldots,p} \frac{  g_{(j)}  }{\sigma \lambda_j}  > 12 \sqrt{\log(1/\delta) }  \right) \le \delta,$$
	which concludes the proof. 
\end{Proof}

\section {Proof of Theorem \ref{cone-condition}}  \label{sec: appendix_cone-condition}
We use the minimality of  $\hat{\B{\beta}}$ and Lemma \ref{lemma-Slope} to derive the cone condition. 

\begin{Proof}
	We assume without loss of generality that $|h_1| \ge \ldots \ge |h_p|$. We define $S_0 = \left\{1,\ldots,k^*\right\}$ as the set of the $k^*$ highest coefficients of $\mathbf{h} =  \hat{\B{\beta}} - \B{\beta}^*$.\\
	
	$\hat{\B{\beta}}$ is the solution of Problem \eqref{learning} hence:
	\begin{align}\label{inf-equation}
	\begin{split}
	&\frac{1}{n} \sum_{i=1}^n f \left( \langle \mathbf{x_i},  \hat{\B{\beta}} \rangle ;  y_i \right) + \lambda \| \hat{\B{\beta}} \|_1  \le \frac{1}{n} \sum_{i=1}^n f \left( \langle \mathbf{x_i},  \B{\beta}^* \rangle ;  y_i \right) + \lambda \| \B{\beta}^* \|_1.
	\end{split}
	\end{align}
	Using the definition of $ \Delta \left(\B{\beta}^*, \mathbf{h} \right) $ as introduced in Theorem \ref{hoeffding-sup}, Equation \eqref{inf-equation} can be written in a more compact form as:
	$$\Delta \left(\B{\beta}^*, \mathbf{h} \right) \le \lambda \| \B{\beta}^* \|_1 - \lambda \| \hat{\B{\beta}} \|_1.$$
	
	Introducing the support $S^*$ of $\B{\beta}^*$ we have
	\begin{align}\label{inf-equation-sup}
	\begin{split}
	\Delta \left(\B{\beta}^*, \mathbf{h} \right)
	&\le  \lambda \| \B{\beta}^*_{S^*} \|_1 - \lambda \| \hat{\B{\beta}}_{S^*} \|_1 -\lambda \| \hat{\B{\beta}}_{(S^*)^c} \|_1\\
	&\le  \lambda \| \mathbf{h}_{S^*} \|_1 -\lambda \| \mathbf{h}_{(S^*)^c}  \|_1\\
	&\le  \lambda \| \mathbf{h}_{S_0} \|_1 -\lambda \| \mathbf{h}_{(S_0)^c}  \|_1,
	\end{split}
	\end{align}
	where this last relation holds by definition of $S_0$. We now want to lower bound $\Delta \left(\B{\beta}^*, \mathbf{h} \right)$. Exploiting the existence of a bounded sub-Gradient $\partial f $ we obtain
	$$\Delta \left(\B{\beta}^*, \mathbf{h} \right)  \ge S \left(\B{\beta}^*, \mathbf{h} \right) :=  \frac{1}{n} \sum_{i=1}^n \partial f \left( \langle \mathbf{x_i},  \B{\beta}^*  \rangle ;  y_i \right) \langle \mathbf{x_i},  \mathbf{h}  \rangle. $$
	
	In addition we have:
	\begin{align*}
	\begin{split}
	|  S \left(\B{\beta}^*, \mathbf{h} \right) | 
	&=\left|  \frac{1}{n} \sum_{i=1}^n \sum_{j=1}^{p} \partial f \left( \langle \mathbf{x_i},  \B{\beta}^*  \rangle ;  y_i \right)   x_{ij} h_j   \right| \\
	&\le \frac{1}{\sqrt{n}}  \sum_{j=1}^{p} \left(   \frac{1}{\sqrt{n} } \left| \sum_{i=1}^n \partial f \left( \langle \mathbf{x_i},  \B{\beta}^*  \rangle ;  y_i \right)  x_{ij} \right|   \right) |h_j|.
	\end{split}
	\end{align*}
	
	Let us define the independent random variables $ g_j =  \frac{1}{\sqrt{n} } \sum_{i=1}^n \partial f \left( \langle \mathbf{x_i},  \B{\beta}^*  \rangle ;  y_i \right)  x_{ij} , \ j=1,\ldots,p.$
	
	Assumption \ref{asuSG} guarantees that $g_1,\ldots,g_p$ are sub-Gaussian with variance $L^2 M^2$. A first upper-bound of the quantity $| S(\mathbf{h}) | $ could be obtained by considering the maximum of the sequence $\left\{g_j\right\}$. However Lemma \ref{upper-bound-sup} gives us a stronger result.\\
	
	Indeed, since $\delta \le 1$ we introduce a non-increasing rearrangement $(g_{(1)}, \ldots, g_{(p)})$  of $(|g_1|, \ldots, |g_p|)$. We recall that $S_0 = \left\{1,\ldots, k^*\right\}$ denotes the subset of indexes of the $k^*$ highest elements of $\mathbf{h}$ and we  use Lemma \ref{upper-bound-sup} to get, with probability at least $1-\frac{\delta}{2}$:
	\begin{align}\label{upper-bound-SG}
	\begin{split}
	|  S \left(\B{\beta}^*, \mathbf{h} \right) | &\le \frac{1}{\sqrt{n}} \sum_{j=1}^p g_{j} |h_{j} | = \frac{1}{\sqrt{n}} \sum_{j=1}^p g_{(j)} |h_{(j)} | = \frac{1}{\sqrt{n}}  \sum_{j=1}^p \frac{ g_{(j)} }{LM \lambda_j} LM \lambda_j | h_{(j)} | \\
	&\le \frac{1}{\sqrt{n}} \sup \limits_{j=1,\ldots,p} \left\{ \frac{ g_{(j)} }{LM \lambda_j} \right\} 
	\sum_{j=1}^p  LM \lambda_j | h_{(j)} | \\
	&\le 12 L M \sqrt{ \frac{\log(2/ \delta)}{n} } \sum_{j=1}^p \lambda_j | h_{(j)} | \textnormal{ with Lemma } \ref{upper-bound-sup}\\
	&\le 12 L M \sqrt{ \frac{\log(2/ \delta)}{n} } \sum_{j=1}^p \lambda_j | h_{j} | \textnormal{ since } \lambda_1\ge \ldots\ge \lambda_p \textnormal{ and } |h_1|\ge \ldots\ge |h_p| \\
	&\le 12 L M \sqrt{ \frac{\log(2/ \delta)}{n} } \left( \sum_{j=1}^{k^*} \lambda_j | h_j |  +  \lambda_{k^*} \sum_{j=k^*}^p | h_j |  \right) \\
	&= 12 L M \sqrt{ \frac{\log(2/ \delta)}{n} } \left( \sum_{j=1}^{k^*} \lambda_j | h_j |  +  \lambda_{k^*} \| \mathbf{h}_{(S_0)^c}  \|_1  \right).
	\end{split}
	\end{align}
	Cauchy-Schwartz inequality leads to:
	\begin{align*}
	\sum_{j=1}^{k^*} \lambda_j | h_j |  &\le \sqrt{\sum_{j=1}^{k^*} \lambda_j ^2 } \| \mathbf{h}_{S_0}  \|_2 \le \sqrt{k^*\log(2pe /k^*)}  \| \mathbf{h}_{S_0}  \|_2,
	\end{align*}
	where we have used the Stirling formula to obtain
	\begin{align*}
	\sum_{j=1}^{k^*} \lambda_j ^2 = \sum_{j=1}^{k^*}  \log(2p/j) &= k^*  \log(2p) - \log(k^* !) \\
	&\le k^* \log(2p) - k^*\log(k^*/ e) = k^* \log(2pe /k^*).
	\end{align*}
	In the statement of Theorem \ref{cone-condition} we have defined $\lambda = 12 \alpha L M \sqrt{ n^{-1} \log(2pe/k^*) \log(2/ \delta)}$. 
	
	Because $\lambda_{k^*} \le  \sqrt{\log(2pe/k^*)}$, Equation \eqref{upper-bound-SG} leads to:
	$$\left|  S \left(\B{\beta}^*, \mathbf{h} \right) \right| \le \frac{1}{\alpha}\lambda \left( \sqrt{k^*}  \| \mathbf{h}_{S_0}  \|_2  +  \| \mathbf{h}_{(S_0)^c}  \|_1 \right)$$
	
	Combined with Equation \eqref{inf-equation-sup}, it holds with probability at least $1-\frac{\delta}{2}:$
	$$ - \frac{\lambda}{\alpha} \left( \sqrt{k^*}  \| \mathbf{h}_{S_0}  \|_2  +  \| \mathbf{h}_{(S_0)^c}  \|_1 \right) \le  \lambda \| \mathbf{h}_{S_0} \|_1 -\lambda \| \mathbf{h}_{(S_0)^c}  \|_1,$$
	which immediately leads to:
	$$ \| \mathbf{h}_{(S_0)^c}  \|_1 \le \frac{\alpha}{\alpha -1}  \| \mathbf{h}_{S_0}  \|_1 + \frac{ \sqrt{k^*}}{\alpha -1}  \| \mathbf{h}_{S_0}  \|_2.$$
	We conclude that $\mathbf{h} \in \Lambda \left(S_0, \ \frac{\alpha}{\alpha -1}, \  \frac{\sqrt{k^*}}{\alpha -1} \right)$ with probability at least $1-\frac{\delta}{2}$.
\end{Proof}

\section {Proof of Theorem \ref{hoeffding-sup}: }  \label{sec: hoeffding-sup}
\begin{Proof}
	Let $k\in \left\{1, \ldots, p\right\}$ and  $S_1, \ldots S_q$ be a partition of $\left\{1, \ldots, p\right\}$ such that $q = \lceil p/k \rceil$ and $|S_j| \le k, \forall j$. We divide the proof of the theorem in 3 steps. We first upper-bound the inner supremum for any sequence of $k$ sparse vectors $\mathbf{z}_{S_1}, \ldots, \mathbf{z}_{S_q}$. We then extend this bound for the supremum over a compact set of sequences through an $\epsilon$-net argument.
	
	\paragraph{Step 1: } Let us fix a sequence  $\mathbf{z}_{S_1}, \ldots, \mathbf{z}_{S_q} \in \mathbb{R}^{p}: \ \Supp(\mathbf{z}_{S_j})  \subset S_j, \forall j$ and $\| \mathbf{z}_{S_j}  \|_1 \le 3R,  \forall j$. 
	
	In particular, $\| \mathbf{z}_{S_j} \|_0 \le k, \forall j$. In the rest of the proof, we define  $\mathbf{z}_{S_0} = \B{0}$ and:
	\begin{equation}  \label{w_l_def}
	\mathbf{w}_{\ell} = \B{\beta}^*  + \sum_{j=1}^{\ell} \mathbf{z}_{S_j}, \forall \ell,
	\end{equation}
	In addition, we  introduce $Z_{i \ell}, \ \forall i,  \ell$ as follows
	$$
	Z_{i \ell} = f \left( \langle \mathbf{x}_i, \mathbf{w}_{\ell}   \rangle ;  y_i \right) -  f \left( \langle \mathbf{x_i},  \mathbf{w}_{\ell-1}  \rangle ;  y_i \right) 
	= f \left( \langle \mathbf{x_i},  \mathbf{w}_{\ell-1} + \mathbf{z}_{S_\ell}  \rangle ;  y_i \right) - f \left( \langle \mathbf{x}_i, \mathbf{w}_{\ell-1}   \rangle ;  y_i \right).  
	$$
	In particular, let us note that:
	\begin{align}
	\begin{split}
	\Delta \left( \mathbf{w}_{\ell-1} , \ \mathbf{z}_{S_\ell} \right) 
	&= \frac{1}{n} \sum_{i=1}^n   f \left( \langle \mathbf{x}_i,  \mathbf{w}_{\ell-1} + \mathbf{z}_{S_{\ell}} \rangle ;  y_i \right)  - \frac{1}{n} \sum_{i=1}^n   f \left( \langle \mathbf{x}_i,   \mathbf{w}_{\ell-1}\rangle ;  y_i \right)  \\
	&= \frac{1}{n} \sum_{i=1}^n \left\{  f \left( \langle \mathbf{x}_i,  \mathbf{w}_{\ell-1} + \mathbf{z}_{S_{\ell}} \rangle ;  y_i \right)  -  f \left( \langle \mathbf{x}_i,   \mathbf{w}_{\ell-1}\rangle ;  y_i \right) \right\}  \\
	&= \frac{1}{n} \sum_{i=1}^n  Z_{i\ell}.
	\end{split}
	\end{align}
	
	Assumption \ref{asu1} guarantees that $f(.,y)$ is L-Lipschitz $\forall y$ then:
	$$|Z_{ij} | \le L \left| \langle \mathbf{x}_i, \mathbf{z}_{S_{\ell}}  \rangle \right|.$$
	Then using  Assumption \ref{asu4}$.1(k)$ on the $k$ sparse vector $\mathbf{z}_{S_{\ell}}$ it holds:
	$$ \left| \Delta \left( \mathbf{w}_{\ell-1} , \ \mathbf{z}_{S_\ell} \right)  \right|
	\le \frac{1}{n} \sum_{i=1}^n  | Z_{i\ell} |
	\le \frac{1}{n} \sum_{i=1}^n L \left| \langle \mathbf{x}_i, \mathbf{z}_{S_{\ell}}  \rangle \right|
	= \frac{L}{n} \| \mathbb{X}  \mathbf{z}_{S_{\ell}} \|_1
	\le \frac{L \mu(k)}{\sqrt{nk} } \| \mathbf{z}_{S_{\ell}} \|_1.
	$$
	
	Hence, with Hoeffding's lemma, the centered bounded random variable $\Delta \left( \mathbf{w}_{\ell-1} , \ \mathbf{z}_{S_\ell} \right) - \mathbb{E}\left( \Delta\left( \mathbf{w}_{\ell-1}  , \ \mathbf{z}_{S_\ell} \right)  \right) $ is sub-Gaussian with variance $\frac{L^2 \mu(k)^2}{nk} \| \mathbf{z}_{S_{\ell}} \|_1^2$. It then hold, $\forall t>0$,
	\begin{align} \label{upper-bound-2k-sparse}
	\begin{split}
	\mathbb{P}\left( \left| \Delta\left(  \mathbf{w}_{\ell-1} , \ \mathbf{z}_{S_\ell} \right)  - \mathbb{E}\left( \Delta\left( \mathbf{w}_{\ell-1}  , \ \mathbf{z}_{S_\ell} \right)  \right) \right|  \ge t  \| \mathbf{z}_{S_\ell} \|_1  \right) 
	&\le 2\exp \left(- \frac{k n t^2 }{ 2 L^2 \mu(k)^2 }  \right).
	\end{split}
	\end{align}
	Equation \eqref{upper-bound-2k-sparse} holds for all values of $\ell$. Thus, an union bound immediately gives:
	\begin{equation} \label{union-upper-bound-2k-sparse}
	\mathbb{P}\left( \sup \limits_{\ell=1,\ldots,q} \left\{   \left| \Delta\left( \mathbf{w}_{\ell-1} , \ \mathbf{z}_{S_\ell} \right)  - \mathbb{E}\left( \Delta\left(  \mathbf{w}_{\ell-1}  , \ \mathbf{z}_{S_\ell} \right)  \right) \right|  - t   \| \mathbf{z}_{S_\ell} \|_1 \right\} \ge 0 \right) \le 2 \left\lceil \frac{p}{k} \right\rceil \exp \left(- \frac{k n t^2 }{ 2 L^2 \mu(k)^2 }  \right).
	\end{equation}

	\paragraph{Step 2: } We  extend the result to any sequence of vectors $\mathbf{z}_{S_1}, \ldots, \mathbf{z}_{S_q} \in \mathbb{R}^{p}: \ \Supp(\mathbf{z}_{S_j})  \subset S_j, \forall j$ and $\| \mathbf{z}_{S_j}  \|_1 \le 3R,  \forall j$ throught an $\epsilon$-net argument.
	\\
	\newline
	We recall that an $\epsilon$-net of a set $\mathcal{I}$ is a subset $\mathcal{N}$ of $\mathcal{I}$ such that each element of $I$ is at a distance at most  $\epsilon$ of $\mathcal{N}$. We know from Lemma 1.18 from \cite{lecture-notes}, that for any $\epsilon \in (0,1)$, the ball $\left\{ \mathbf{z} \in \mathbb{R}^d: \ \|  \mathbf{z} \|_1 \le R  \right\}$ has an $\epsilon$-net of cardinality $| \mathcal{N} | \le \left(\frac{2R+1}{\epsilon} \right)^d$ -- the $\epsilon$-net is defined in term of L1 norm. In addition, by following the proof of the lemma, we can create this set such that it contains $\B{0}$.
	\smallskip
	
	Consequently, we use Equation \eqref{union-upper-bound-2k-sparse} on a product of $\epsilon$-nets $\mathcal{N}_{k,R} = \prod \limits_{\ell=1}^q \mathcal{N}_{k,R}^{\ell}$.  Each $\mathcal{N}_{k,R}^{\ell}$ is an $\epsilon$-net of the bounded sets of $k$ sparse vectors $ \mathcal{I}_{k,R}^{\ell} = \left\{ \mathbf{z}_{S_{\ell}} \in \mathbb{R}^{p}: \ \Supp(\mathbf{z}_{S_{\ell}})  \subset S_{\ell} \ ; \ \| \mathbf{z}_{S_{\ell}}  \|_1 \le 3R\right\}$ which contains $\B{0}_{S_{\ell}}$. We note $\mathcal{I}_{k,R} = \prod \limits_{\ell=1}^q \mathcal{I}_{k,R}^{\ell}$. It then holds:
	\begin{align} \label{sup-epsilon-net}
	\begin{split}
	&\mathbb{P}\left( \sup \limits_{ \left( \mathbf{z}_{S_1}, \ldots, \mathbf{z}_{S_q}\right) \in \mathcal{N}_{k,R}  } \ \ \left\{ \sup \limits_{\ell=1,\ldots,q} \left\{   \left| \Delta\left( \mathbf{w}_{\ell-1} , \ \mathbf{z}_{S_\ell} \right)  - \mathbb{E}\left( \Delta\left(  \mathbf{w}_{\ell-1}  , \ \mathbf{z}_{S_\ell} \right)  \right) \right|  - t  \| \mathbf{z}_{S_\ell} \|_1 \right\} \ge 0 \right\} \right)\\
	&\le 2 \left\lceil \frac{p}{k} \right\rceil  \left(\frac{6R+1}{\epsilon} \right)^{k} \left\lceil \frac{p}{k} \right\rceil  \exp \left(- \frac{k n t^2 }{ 2 L^2 \mu(k)^2 } \right) \le 2 \left( \frac{2p}{k} \right)^2 \left(\frac{6R+1}{\epsilon} \right)^{k} \exp \left(- \frac{k n t^2 }{ 2 L^2 \mu(k)^2 } \right).
	\end{split}
	\end{align}
	
	\paragraph{Step 3: } We now extend Equation \eqref{sup-epsilon-net} to control any vector in $\mathcal{I}_{k,R}$.  For $\mathbf{z}_{S_1}, \ldots, \mathbf{z}_{S_q} \in \mathcal{I}_{k,R}$, there exists $\tilde{\mathbf{z}}_{S_1}, \ldots, \tilde{\mathbf{z}}_{S_q}\in \mathcal{N}_{k,R}$ such that $\| \mathbf{z}_{S_\ell} - \tilde{\mathbf{z}}_{S_\ell} \|_1 \le \epsilon, \forall \ell.$ Similarly to Equation \eqref{w_l_def}, we define:
	$$ \tilde{\mathbf{w}}_{\ell} = \B{\beta}^*  + \sum_{j=1}^{\ell} \tilde{\mathbf{z}}_{S_j}, \forall \ell.$$ 
	For a given $t$, let us define 
	$$f_t \left( \mathbf{w}_{\ell-1} , \ \mathbf{z}_{S_{\ell}} \right) = \left| \Delta \left( \mathbf{w}_{\ell-1} , \ \mathbf{z}_{S_{\ell}} \right) - \mathbb{E} \left( \mathbf{w}_{\ell-1} , \ \mathbf{z}_{S_{\ell}} \right) \right| - t  \| \mathbf{z}_{S_{\ell}} \|_1, \forall \ell.$$
	
	We fix $\ell_0(t)$ such that $\ell_0 \in \argmax \limits_{\ell=1,\ldots,q}  \left\{   f_{7t} \left( \mathbf{w}_{\ell-1} , \ \mathbf{z}_{S_{\ell}} \right)  \right\}$. The choice of $7t$ will be justified later. We fix $t$ and will just note $\ell_0=\ell_0(t)$ when no confusion can be made. 
	\\
	\newline
	With Assumption \ref{asu1} we obtain:
	\begin{align} \label{help-step3}
	\begin{split}
	&\left| \Delta \left( \mathbf{w}_{\ell_0-1} , \ \mathbf{z}_{S_{\ell_0}} \right)
	- \Delta\left( \tilde{\mathbf{w}}_{\ell_0-1} , \ \tilde{\mathbf{z}}_{S_{\ell_0}} \right) \right| \\
	&= \frac{1}{n} \left| \sum_{i=1}^n f \left( \langle \mathbf{x}_i,  \mathbf{w}_{\ell_0}   \rangle ;  y_i \right) -\sum_{i=1}^n f \left( \langle \mathbf{x}_i,  \tilde{\mathbf{w}}_{\ell_0}  \rangle ;  y_i \right) +  \sum_{i=1}^n  f \left( \langle \mathbf{x}_i,  \tilde{\mathbf{w}}_{\ell_0-1}  \rangle ;  y_i \right)    - \sum_{i=1}^n  f \left( \langle \mathbf{x}_i,  \mathbf{w}_{\ell_0-1}  \rangle ;  y_i \right) \right| \\
	&\le \frac{1}{n} \sum_{i=1}^n L \left| \langle \mathbf{x}_i, \mathbf{w}_{\ell_0} - \tilde{\mathbf{w}}_{\ell_0} \rangle  \right| + \frac{1}{n} \sum_{i=1}^n L \left| \langle \mathbf{x}_i, \mathbf{w}_{\ell_0-1} - \tilde{\mathbf{w}}_{\ell_0-1} \rangle  \right|  \\
	&= \frac{1}{n} \sum_{i=1}^n L \left|  \sum_{\ell=1}^{\ell_0}  \langle \mathbf{x}_i, \mathbf{z}_{S_\ell} - \tilde{\mathbf{z}}_{S_\ell} \rangle  \right| + \frac{1}{n} \sum_{i=1}^n L \left| \sum_{\ell=1}^{\ell_0-1}  \langle \mathbf{x}_i,\mathbf{z}_{S_\ell} - \tilde{\mathbf{z}}_{S_\ell} \rangle  \right|  \\
	&\le \frac{2}{n} \sum_{i=1}^n \sum_{\ell=1}^{q} L \left| \langle \mathbf{x}_i, \mathbf{z}_{S_\ell} - \tilde{\mathbf{z}}_{S_\ell} \rangle  \right| \\
	&= \frac{2}{\sqrt{n} } \sum_{\ell=1}^{q}   \frac{L}{\sqrt{n} }  \left\| \mathbf{X} \right( \mathbf{z}_{S_\ell} - \tilde{\mathbf{z}}_{S_\ell} \left) \right\|_1 \\
	&\le  \frac{2}{ \sqrt{n} } \sum_{\ell=1}^{q}  \frac{L}{\sqrt{k} }  \mu(k) \ \left\| \mathbf{z}_{S_\ell} - \tilde{\mathbf{z}}_{S_\ell} \right\|_1 \\
	&\le  \frac{2p}{ k\sqrt{k n} }  L \mu(k) \epsilon \le  \eta \epsilon.
	\end{split}
	\end{align}
	where $\eta = \frac{2L \mu(k)}{ \sqrt{n} }$ and we have used Assumption \ref{asu5}.1($k$). It then holds:
	\begin{align*} 
	\begin{split}
	f_t \left( \tilde{\mathbf{w}}_{\ell_0-1} , \ \tilde{\mathbf{z}}_{S_{\ell_0}} \right)  
	&\ge f_t \left( \mathbf{w}_{\ell_0-1} , \ \mathbf{z}_{S_{\ell_0}} \right) -  \left| \Delta \left( \mathbf{w}_{\ell_0-1} , \ \mathbf{z}_{S_{\ell_0}} \right) -  \Delta\left( \tilde{\mathbf{w}}_{\ell_0-1} , \ \tilde{\mathbf{z}}_{S_{\ell_0}} \right)  \right|  \\
	&\ \ \ - \left| \mathbb{E}\left(  \Delta \left( \mathbf{w}_{\ell_0-1} , \ \mathbf{z}_{S_{\ell_0}}\right) -  \Delta\left( \tilde{\mathbf{w}}_{\ell_0-1} , \ \tilde{\mathbf{z}}_{S_{\ell_0}} \right) \right)\right|   - t  \| \mathbf{z}_{S_{\ell_0}} - \tilde{\mathbf{z}}_{S_{\ell_0}}  \|_1  \\
	&\ge f_t \left( \mathbf{w}_{\ell_0-1} , \ \mathbf{z}_{S_{\ell_0}} \right)  -  2 \eta \epsilon - t \epsilon.
	\end{split}
	\end{align*}
	
	\textbf{Case 1:} Let us assume that $\| \mathbf{z}_{S_{\ell_0}} \|_1 \ge \epsilon / 2$ and that $t\ge \eta$, then we have:
	\begin{equation}\label{step3-eps}
	f_t \left( \tilde{\mathbf{w}}_{\ell_0-1} , \ \tilde{\mathbf{z}}_{S_{\ell_0}} \right)  \ge f_t \left( \mathbf{w}_{\ell_0-1} , \ \mathbf{z}_{S_{\ell_0}} \right)  -  2(2\eta + t ) \| \tilde{\mathbf{z}}_{S_{\ell_0}} \|_1 \ge f_{7t} \left( \mathbf{w}_{\ell_0-1} , \ \mathbf{z}_{S_{\ell_0}} \right). 
	\end{equation}
	
	\textbf{Case 2:} We now assume $\| \mathbf{z}_{S_{\ell_0}} \|_1 \le \epsilon / 2$. Since $\B{0}_{S_{\ell_0} } \in  \mathcal{N}_{k,R}$ we derive similarly to Equation  \eqref{help-step3}:
	$$\left| \Delta \left( \mathbf{w}_{\ell_0-1} , \ \mathbf{z}_{S_{\ell_0}} \right)
	- \Delta\left( \mathbf{w}_{\ell_0-1} , \ \mathbf{0}_{S_{\ell_0}} \right) \right| 
	\le \frac{L \mu(k) }{\sqrt{n k} } \left\| \mathbf{z}_{S_{\ell_0} } \right\|_1,$$
	which then implies that:
	$$f_{7t} \left( \mathbf{w}_{\ell_0-1} , \ \mathbf{z}_{S_{\ell_0}} \right) 
	\le f_{7t} \left( \mathbf{w}_{\ell_0-1} , \ \mathbf{0}_{S_{\ell_0}}  \right) + \frac{2 L \mu(k) }{\sqrt{n k}} \left\| \mathbf{z}_{S_{\ell_0} } \right\|_1 - 7t  \left\| \mathbf{z}_{S_{\ell_0} } \right\|_1,$$
	and this quantity is smaller than $f_{7t} \left( \mathbf{w}_{\ell_0-1} , \ \mathbf{0}_{S_{\ell_0}}  \right)$ as long as $7 t \ge \frac{2 L \mu(k) }{\sqrt{n k}}$. The latter condition is satisfied if $t\ge \eta$.
	
	In this case, we can define a new $\tilde{\ell}_0$ for the sequence $\mathbf{z}_{S_1}, \ldots, \mathbf{z}_{S_{\ell_0-1}}, \mathbf{0}_{S_{\ell_0}}, \mathbf{z}_{S_{\ell_0+1}}, \ldots, \mathbf{z}_{S_q}$.  After a finite number of iteration, by using the result in Equation \eqref{step3-eps} and the definition of $\ell_0$, we finally get that  $f_{7t} \left( \mathbf{w}_{\ell_0-1} , \ \mathbf{z}_{S_{\ell_0}} \right) \le f_t \left( \tilde{\mathbf{w}}_{\ell_0-1} , \ \tilde{\mathbf{z}}_{S_{\ell_0}} \right)$ for some $\tilde{\mathbf{z}}_{S_1}, \ldots, \tilde{\mathbf{z}}_{S_q}\in \mathcal{N}_{k,R}$.
	\\
	\newline
	As a consequence of cases 1 and 2, we obtain: $\forall t \ge \eta, \ \forall \mathbf{z}_{S_1}, \ldots, \mathbf{z}_{S_q} \in \mathcal{I}_{k,R}, \ \exists \tilde{\mathbf{z}}_{S_1}, \ldots, \tilde{\mathbf{z}}_{S_q}\in \mathcal{N}_{k,R}$:
	$$\sup \limits_{\ell=1,\ldots,q}  f_{7t} \left( \mathbf{w}_{\ell-1} , \ \mathbf{z}_{S_{\ell}} \right) = f_{7t} \left( \mathbf{w}_{\ell_0-1} , \ \mathbf{z}_{S_{\ell_0}} \right)  \le f_{t} \left( \tilde{\mathbf{w}}_{\ell_0-1} , \ \tilde{\mathbf{z}}_{S_{\ell_0}} \right) \le \sup \limits_{\ell=1,\ldots,q} f_{t} \left( \tilde{\mathbf{w}}_{\ell-1} , \ \tilde{\mathbf{z}}_{S_{\ell}} \right).$$
	
	This last relation is equivalent to saying that $\forall t \ge 7\eta$:
	\begin{equation}\label{supremum-domination}
	\sup \limits_{ \mathbf{z}_{S_1}, \ldots, \mathbf{z}_{S_q} \in \mathcal{I}_{k,R} } \left\{    \sup \limits_{\ell=1,\ldots,q}  f_t \left( \mathbf{w}_{\ell-1} , \ \mathbf{z}_{S_{\ell}} \right)  \right\}  \le 
	\sup \limits_{ \mathbf{z}_{S_1}, \ldots, \mathbf{z}_{S_q} \in \mathcal{N}_{k,R} } \left\{    \sup \limits_{\ell=1,\ldots,q}  f_{t/7} \left( \tilde{\mathbf{w}}_{\ell-1} , \ \tilde{\mathbf{z}}_{S_{\ell}}, \right) \right\}.
	\end{equation}

	As a consequence, we have $\forall t\ge 7\eta$:
	\begin{align}  \label{sup-compact-set}.
	\begin{split}
	&\mathbb{P}\left( \sup \limits_{ \mathbf{z}_{S_1}, \ldots, \mathbf{z}_{S_q} \in \mathcal{I}_{k,R} } \ \ \sup \limits_{\ell=1,\ldots,q} \left\{   \left| \Delta\left( \mathbf{w}_{\ell-1} , \ \mathbf{z}_{S_\ell} \right)  - \mathbb{E}\left( \Delta\left(  \mathbf{w}_{\ell-1}  , \ \mathbf{z}_{S_\ell} \right)  \right) \right|  - t  \| \mathbf{z}_{S_\ell} \|_1 \right\} \ge 0  \right)\\
	&\le \mathbb{P}\left( \sup \limits_{ \mathbf{z}_{S_1}, \ldots, \mathbf{z}_{S_q} \in \mathcal{N}_{k,R} } \ \ \sup \limits_{\ell=1,\ldots,q} \left\{   \left| \Delta\left( \mathbf{w}_{\ell-1} , \ \mathbf{z}_{S_\ell} \right)  - \mathbb{E}\left( \Delta\left(  \mathbf{w}_{\ell-1}  , \ \mathbf{z}_{S_\ell} \right)  \right) \right|  - \frac{t}{7}  \| \mathbf{z}_{S_\ell} \|_1 \right\} \ge 0  \right)\\
	&\le 2 \left( \frac{2p}{k} \right)^2 \left(\frac{6R+1}{\epsilon} \right)^{k} \exp \left(- \frac{ k n (t/7)^2 }{ 2 L^2 \mu(k)^2 } \right)\\
	& \le \left( \frac{4p}{k} \right)^2 3^{k}  \exp \left(- \frac{k n t^2 }{98 L^2 \mu(k)^2 } \right) \textnormal{ by fixing } \epsilon = 2R \textnormal{ and since } R \ge 1.
	\end{split}
	\end{align}
	Thus we select $t$ such that $t\ge 7\eta$ and that the condition $t^2 \ge \frac{98 L^2 \mu(k)^2 }{2kn}\left[k \log(3) + 2 \log\left( \frac{4 p}{k} \right) + \log\left( \frac{2}{\delta}\right) \right]$ holds \footnote{A somewhat faster proof would have consisted in fixing $\epsilon = 2R$ in the definition of the $\epsilon$-net -- of size now bounded by $3^k$ -- and in noting that because of the L1-constraint, each element $\mathbf{z}_{S_{\ell}}$ is at a distance at most $R = \| \mathbf{z}_{S_{\ell}} \|_1 / 2$ of its closest neighborhood in the $\epsilon$-net. However we prefer the more general proof presented.}. To this end, we define: 
	$$\tau = 14 L  \mu(k) \sqrt{\frac{\log(3) }{n} + \frac{\log\left( 4 p/k \right) }{nk}   +  \frac{ \log\left( 2/\delta\right)  }{nk}  } \ge 7 \eta.$$
	We conclude that with probability at least $1-\frac{\delta}{2}$:
	$$ \sup \limits_{ \mathbf{z}_{S_1}, \ldots, \mathbf{z}_{S_q} \in \mathcal{I}_{k,R} } \left\{\sup \limits_{\ell=1,\ldots,q} \left\{   \left| \Delta\left( \mathbf{w}_{\ell-1} , \ \mathbf{z}_{S_\ell} \right)  - \mathbb{E}\left( \Delta\left(  \mathbf{w}_{\ell-1}  , \ \mathbf{z}_{S_\ell} \right)  \right) \right|  - \tau  \left(\| \mathbf{z}_{S_\ell} \|_1\vee \eta\right) \right\}  \right\}  \le 0.$$
\end{Proof}

\section {Proof of Theorem \ref{restricted-strong-convexity}: }  \label{sec: appendix_restricted-strong-convexity}
\begin{Proof} 
	The proof is divided in two steps. First, we lower-bound the quantity $\Delta\left(\B{\beta}^*, \mathbf{h} \right) $ by using a decomposition of $\left\{1, \ldots, p\right\}$  and applying Theorem \ref{hoeffding-sup}. Second, we consider the cone condition derived in Theorem \ref{cone-condition} and use the restricted eigenvalue condition presented in Assumption \ref{asu4}$.2$.

	\paragraph{Step 1:} 
	Let us fix the partition $S_1= \left\{1, \ldots, k^*\right\}, S_2= \left\{k^*+1, \ldots, 2k^*\right\},\ldots, S_q$ of $\left\{1, \ldots, p\right\}$  -- with $q = \lceil p/k^* \rceil$. Thus it holds $|S_j| \le k^*, \forall j$ and we can use Theorem \ref{hoeffding-sup}. We define the corresponding sequence $\mathbf{h}_{S_1}, \ldots, \mathbf{h}_{S_q}$ of $k^*$ sparse vectors corresponding to the decomposition of $\mathbf{h}=\hat{\B{\beta}}-\B{\beta}^*$ on the partition and note that:
	\begin{align} \label{decomposition}
	\begin{split}
	\Delta(\B{\beta}^*, \mathbf{h}) &=\frac{1}{n} \sum_{i=1}^n f \left( \langle \mathbf{x_i},  \B{\beta}^* + \mathbf{h}  \rangle ;  y_i \right)  - \frac{1}{n} \sum_{i=1}^n  f \left( \langle \mathbf{x_i},  \B{\beta}^* \rangle ;  y_i \right) \\
	&=\frac{1}{n} \sum_{i=1}^n f \left( \langle \mathbf{x_i},  \B{\beta}^* + \sum_{j=1}^{q} \mathbf{h}_{S_j}   \rangle ;  y_i \right)  - \frac{1}{n} \sum_{i=1}^n  f \left( \langle \mathbf{x_i},  \B{\beta}^* \rangle ;  y_i \right)\\
	&=  \sum_{\ell=1}^{q} \left\{ \frac{1}{n} \sum_{i=1}^n f \left( \langle \mathbf{x_i},  \B{\beta}^* + \sum_{j=1}^{\ell} \mathbf{h}_{S_j}   \rangle ;  y_i \right)  - \frac{1}{n} \sum_{i=1}^n  f \left( \langle \mathbf{x_i},  \B{\beta}^*  + \sum_{j=1}^{\ell-1} \mathbf{h}_{S_j}  \rangle ;  y_i \right)  \right\}\\
	&= \sum_{\ell=1}^{q} \Delta \left( \B{\beta}^*  + \sum_{j=1}^{\ell-1} \mathbf{h}_{S_j} , \ \mathbf{h}_{S_\ell} \right)\\
	&= \sum_{\ell=1}^{q} \Delta \left( \mathbf{w}_{\ell-1} , \ \mathbf{h}_{S_\ell} \right).	
	\end{split}
	\end{align}
	where we have defined  $\mathbf{w}_{\ell} = \B{\beta}^*  + \sum_{j=1}^{\ell} \mathbf{z}_{S_j}, \forall \ell$ and $\mathbf{z}_{S_0}=\B{0}$ as in the proof of Theorem \ref{hoeffding-sup}. Consequently, since $ \| \mathbf{h}_{S_\ell} \|_0 \le k^*$ and $\| \mathbf{h}_{S_\ell} \|_1 \ge R, \ \forall \ell$,  it holds  with probability at least $1 - \frac{\delta}{2}$:
	$$\left| \Delta \left( \mathbf{w}_{\ell-1} ,  \mathbf{h}_{S_\ell} \right) - \mathbb{E} \left( \mathbf{w}_{\ell-1} ,  \mathbf{h}_{S_\ell}  \right) \right|
	\ge \tau \| \mathbf{h}_{S_\ell}  \|_1, \forall \ell, $$
	where $\tau  = \tau(k^*) = 14 L  \mu(k^*) \sqrt{\frac{\log(3) }{n} + \frac{\log\left( 4 p/k^* \right) }{nk^*}   +  \frac{ \log\left( 2/\delta\right)  }{nk^*}  } $ is fixed in the rest of the proof. 
	\smallskip
	
	As a result, following Equation \eqref{decomposition}, we have:
	\begin{align} \label{lower-bound-expectation-1}
	\begin{split}
	\Delta(\B{\beta}^* , \mathbf{h})
	&\ge   \sum_{\ell=1}^{q}  \left\{ \mathbb{E} \left( \mathbf{w}_{\ell-1}  ,  \mathbf{h}_{S_\ell}  \right) -  \tau \| \mathbf{h}_{S_\ell}  \|_1 \right\}\\
	&= \mathbb{E} \left( \sum_{\ell=1}^{q} \Delta \left( \mathbf{w}_{\ell-1} , \ \mathbf{h}_{S_\ell} \right) \right) - \sum_{\ell=1}^q \tau \| \mathbf{h}_{S_\ell}  \|_1 \\
	&= \mathbb{E} \left( \Delta(\B{\beta}^* , \mathbf{h}) \right) -  \tau \| \mathbf{h}  \|_1.
	\end{split}
	\end{align}
	In addition, we have:
	$$
	\mathbb{E} \left( \Delta(\B{\beta}^* , \mathbf{h}) \right)  = 
	\frac{1}{n} \sum_{i=1}^n \mathbb{E} \left\{  f \left( \langle \mathbf{x_i},  \B{\beta}^* + \mathbf{h}  \rangle ;  y_i \right)  -  f \left( \langle \mathbf{x_i},  \B{\beta}^*  \rangle ;  y_i \right)  \right\} =L(\B{\beta}^* + \mathbf{h}) - L(\B{\beta}^*).
	$$
	Consequently, we conclude that with probability at least $1 - \frac{\delta}{2}$:
	\begin{equation} \label{lower-bound-expectation}
	\Delta(\B{\beta}^* , \mathbf{h})\ge L(\B{\beta}^* + \mathbf{h}) - L(\B{\beta}^*) -   \tau \| \mathbf{h}  \|_1.
	\end{equation}
	
	\bigskip
	
	\paragraph{Step 2:} 
	We now lower-bound the right-hand side of Equation  \eqref{lower-bound-expectation}. Since $\mathcal{L}$ is twice differentiable, a Taylor development around $\B{\beta}^*$ gives:
	$$\mathcal{L}(\B{\beta}^* + \mathbf{h}) - \mathcal{L}(\B{\beta}^*) 
	= \nabla \mathcal{L}(\B{\beta}^*)^T\mathbf{h}  + \frac{1}{2} \mathbf{h}^T \nabla^2 \mathcal{L}(\B{\beta}^*)^T\mathbf{h} + o \left( \| \mathbf{h}\|_2 \right).$$
	The optimality of $\B{\beta}^*$ implies $\nabla L(\B{\beta}^*)=0$. In addition, Theorem \ref{cone-condition} states that $\mathbf{h} \in \Lambda \left(S_0, \gamma_1,  \gamma_2 \right)$ with probability at least $1- \frac{\delta}{2}$. Consequently, we can use the restricted eigenvalue condition defined in Assumption \ref{asu4}$.2(k^*, \gamma)$. However we do not want to keep the term $o \left( \| \mathbf{h} \|_2 \right)$ as it can hide non trivial dependencies. 
	\medskip
	
	\textbf{Case 1:} If $\| \mathbf{h} \|_2 \le r(k^*)$ -- where  $r(k^*, \gamma)$ is shorthanded $r(k^*)$ and is the maximum radius introduced in the growth condition Assumption \ref{asu5}$.2$ -- then by the result of Theorem \ref{cone-condition} and Assumption  \ref{asu4}$.2(k, \gamma)$, it holds with probability at least  $1- \frac{\delta}{2}$:
	\begin{equation}\label{case1-LB}
	\mathcal{L}(\B{\beta}^* + \mathbf{h}) - \mathcal{L}(\B{\beta}^*) \ge \frac{1}{4}  \kappa(k^*) \|\mathbf{h} \|_2^2.
	\end{equation}
	
	\textbf{Case 2:} If now $\| \mathbf{h} \|_2 \ge r(k^*)$, then using the convexity of $\mathcal{L}$ thus of $t \to \mathcal{L}\left( \B{\beta}^* + t \mathbf{h} \right)$, we similarly obtain with the same probability:
	\begin{align} \label{trick}
	\begin{split}
	\mathcal{L}(\B{\beta}^* + \mathbf{h}) - \mathcal{L}(\B{\beta}^*) 
	& \ge  \frac{\|\mathbf{h} \|_2}{r(k^*) } \left\{ \mathcal{L} \left(\B{\beta}^* + \frac{r(k^*) }{\|\mathbf{h} \|_2} \mathbf{h} \right)   - \mathcal{L}	(\B{\beta}^* )  \right\} \textnormal{by convexity} \\
	& \ge  \frac{\|\mathbf{h} \|_2}{r(k^*) } \inf \limits_{ \substack{\mathbf{z}: \  \mathbf{z} \in \Lambda(S_0, \gamma_1, \gamma_2) \\ \| \mathbf{z} \|_2 = r(k^*)}  }  \left\{ \mathcal{L}(\B{\beta}^* + \mathbf{z} )   - \mathcal{L}(\B{\beta}^*)  \right\} \\
	&\ge    \frac{\|\mathbf{h} \|_2}{r(k^*)} \ \frac{1}{4} \kappa(k^*)  r(k^*)^2 = \frac{1}{4} \kappa(k^*) r(k^*)  \|\mathbf{h} \|_2.\\
	\end{split}
	\end{align}
	
	Combining Equations \eqref{lower-bound-expectation}, \eqref{case1-LB} and \eqref{trick}, we conclude that with probability at least $1-\delta$ the following restricted strong convexity holds:
	\begin{equation}\label{proof-rsc}
	\Delta(\mathbf{h}) \ge \frac{1}{4}  \kappa(k^*) \|\mathbf{h} \|_2^2 \wedge   \frac{1}{4} \kappa(k^*) r(k^*)  \|\mathbf{h} \|_2 -   \tau  \| \mathbf{h}\|_1.
	\end{equation}
\end{Proof}

\newpage
\section {Proof of Theorem \ref{main-results} } \label{sec: appendix_main-results}

\begin{Proof} We now prove our main Theorem  \ref{main-results}.
	We recall that $S_0$ has been defined as the subset of the $k^*$ highest elements of $\mathbf{h}$. Following Equation \eqref{inf-equation} and the restricted strong convexity derived in Theorem \ref{restricted-strong-convexity} (Equation \eqref{proof-rsc}), it holds with probability at least $1-\delta$:
	\begin{equation}\label{use-tau-lambda}
	\begin{split}
	&\frac{1}{4}  \kappa(k^*) \left\{ \|\mathbf{h}\|_2^2 \wedge r(k^*) \|\mathbf{h}\|_2 \right\} \\
	&\le \tau \|\mathbf{h}\|_1 + \lambda \| \mathbf{h}_{S^*} \|_1 -\lambda \| \mathbf{h}_{(S^*)^c}  \|_1 \\
	&= \tau \left( \| \mathbf{h}_{S_0} \|_1 + \| \mathbf{h}_{(S_0)^c} \|_1 \right) + \lambda\sqrt{k^*} \| \mathbf{h}_{S_0} \|_2\\
	&\le \tau \left( \| \mathbf{h}_{S_0} \|_1 +  \frac{\alpha}{\alpha -1}  \| \mathbf{h}_{S_0}  \|_1 + \frac{ \sqrt{k^*}}{\alpha -1}  \| \mathbf{h}_{S_0}  \|_2 \right) + \lambda\sqrt{k^*} \| \mathbf{h}_{S_0} \|_2 \text{  since  }\mathbf{h} \in \Lambda \left(S_0, \gamma_1,  \gamma_2 \right) \\
	&\le \frac{2\alpha - 1}{\alpha -1}  \tau \sqrt{k^*}  \|\mathbf{h}_{S_0} \|_2 +  \frac{\tau \sqrt{k^*} }{\alpha -1}   \|\mathbf{h}_{S_0} \|_2 +  \lambda\sqrt{k^*} \| \mathbf{h}_{S_0} \|_2\\
	& \ \ \ \ \  \ \ \ \text{ from Cauchy-Schwartz inequality on the } k^* \text{ sparse vector } \mathbf{h}_{S_0}\\
	&\le \left( \frac{2\alpha   }{\alpha -1}  \tau + \lambda \right) \sqrt{k^*}  \|\mathbf{h} \|_2.
	\end{split}
	\end{equation}
	
	With the definitions of $\tau$ and $\lambda$  as in the Theorems \ref{cone-condition} and  \ref{hoeffding-sup}, Equation \eqref{use-tau-lambda} leads to: 
	\begin{equation*}
	\begin{split}
	\frac{1}{4}  \kappa(k^*) \left\{ \|\mathbf{h}\|_2 \wedge r(k^*)\right\}  
	&\le 12 \alpha L M \sqrt{ \frac{ k^* \log(2pe/k^*) }{n}  \log(2/ \delta) }\\
	&+ \frac{28\alpha }{\alpha -1}  L  \mu(k^*) \sqrt{\frac{\log(3) }{n} + \frac{\log\left( 4 p/k \right) }{nk}   +  \frac{ \log\left( 2/\delta\right)  }{nk}  }.
	\end{split}
	\end{equation*}
	
	Exploiting Assumption \ref{asu5}($ k^*,  \gamma, \delta)$, and using that $\alpha \ge 2$, we obtain with probability at least $1 - \delta$:
	\begin{equation*}
	\begin{split}
	&\|\mathbf{h}\|_2^2 \lesssim 
	\left( \frac{ \alpha L M }{\kappa(k^*)} \right)^2  \frac{ k^* \log\left( p/k^* \right) \log\left( 2/\delta \right) }{n} + \left( \frac{ \alpha L \mu(k^*)  }{ \kappa(k^*)} \right)^2   \frac{ \log(3) + \log\left(p /k^* \right) / k^* + \log\left( 2/ \delta \right) }{n}.
	\end{split}
	\end{equation*}
	which concludes the proof.
\end{Proof}

\section {Proof of Corollary \ref{main-corollary} } \label{sec: appendix_main-corollary}

\begin{Proof}
	In order to derive the bound in expectation, we define the bounded random variable: 
	$$ Z =  \frac{\kappa(k^*)^{2}}{\alpha^2 L^2} \| \hat{\B{\beta}}  - \B{\beta}^*\|_2^2.$$
	Since Assumption \ref{asu5}($k^*,  \gamma, \delta_0$) is satisfied for a small enough $\delta_0$, and by assuming $\log(3) \le k^*$  we can fix $C$ such that $\forall \delta \in \left(0, 1 \right)$, it holds with probability at least $1-\delta$:
	\begin{equation*}
	Z \le C H \left\{ \mu(k^*)^2 + M^2  \log(2/\delta)   \right\} + C \frac{\mu(k^*)^2 }{n} \log(2/\delta) \ \text{ where } \ H= \frac{k^*\log(p/k^*)}{n}.
	\end{equation*}
	
	Then it holds $\forall t \ge t_0 = \log(4):$
	$$\mathbb{P}\left( Z/C \ge H \mu(k^*)^2 + HM^2t +  \frac{\mu(k^*)^2 }{n} t \right) \le 2e^{-t}.$$
	
	Let $q_0 = HM^2t_0 +  \frac{\mu(k^*)^2 }{n} t_0 $, then $\forall q \ge q_0$
	\begin{align*}
	\mathbb{P}\left( Z/C \ge H \mu(k^*)^2 +q \right) &\le 2\exp\left( - \frac{n}{nH M^2 + \mu(k^*)^2} \ q \right) \le 2\exp\left( -  \frac{q}{H M^2} \right).
	\end{align*}
	
	Consequently, by integration we have:
	\begin{align} 
	\begin{split}
	\mathbb{E}(Z) &= \displaystyle \int_0^{+ \infty}  C\mathbb{P}\left( |Z| /C \ge q \right)dq\\
	&= \displaystyle \int_0^{+ \infty}  C\mathbb{P}\left( |Z| /C \ge H \mu(k^*)^2 + q \right)dq + C H \mu(k^*)^2 \\
	&\le \displaystyle \int_{q_0}^{+ \infty}  C \mathbb{P}\left( |Z| /C \ge H \mu(k^*)^2 + q\right) dq + C q_0 + CH\mu(k^*)^2  \\
	&\le \displaystyle \int_{q_0}^{+ \infty}  2C e^{ -\frac{q}{H M^2}   }dq +  C q_0 + CH\mu(k^*)^2\\
	&\le 2CHM^2 e^{ -\frac{q_0}{H M^2}  }  +  C q_0 + C H\mu(k^*)^2 \\
	&\le 2 CHM^2 +  CHM^2 \log(4) +  C\frac{\mu(k^*) }{n}\log(4)  + C H\mu(k^*)^2 \\
	&\le C_1 H ( \mu(k^*)^2 +M^2  ) 
	\end{split}
	\end{align}
	for some universal constant $C_1$, since $H\gg n^{-1} \mu(k^*)$. Hence we conclude:
	$$\mathbb{E}\left( \| \hat{\B{\beta}} - \B{\beta}^*  \|_2^2 \right)  \lesssim   \left( \frac{\alpha L}{\kappa(k^*)} \right)^2( \mu(k^*)^2 +M^2  ) \frac{k^* \log\left(p /k^* \right)}{n} .$$
\end{Proof}

\section {Proof of Theorem \ref{lipschitz} } \label{sec: appendix_lipschitz}

\begin{Proof} 
	We fix $ \tau >0$ and denote $\mathbb{X} = (\mathbf{X}_1, \ldots, \mathbf{X}_p) \in \mathbb{R}^{n \times p}$ the design matrix.
	
	For $\B{\beta} \in \mathbb{R}^{p}$, we define $\textbf{w}^{\tau}(\B{\beta}) \in \mathbb{R}^n$ by:
	$$w^{\tau}_i(\B{\beta}) = \min\left( 1, \frac{1}{2\tau} | z_i | \right) \sign(z_i ), \ \forall i$$
	where $z_i = 1 - y_i \mathbf{x}_i^T\B{\beta}, \ \forall i$. We easily check that
	$$\pmb{w}^{\tau}(\B{\beta}) =   \argmax \limits_{\|w\|_{\infty} \le 1} \frac{1}{2n} \sum \limits_{i=1}^n \left( z_i  + w_i z_i \right) - \frac{\tau}{2n} \|w\|_2^2. $$
	
	Then the gradient of the smooth hinge loss is
	\begin{equation*} \label{smoothed-gradient}
	\nabla g^{\tau}( \B{\beta} ) = - \frac{1}{2n} \sum \limits_{i=1}^{n}(1+w_i^{\tau}(\B{\beta})) y_i \mathbf{x}_i \in \mathbb{R}^{p}.
	\end{equation*}

	For every couple $\B{\beta}, \pmb{\gamma}  \in \mathbb{R}^{p}$ we have:
	\begin{equation} \label{diff-gradient}
	\nabla g^{\tau}(\B{\beta}) - \nabla g^{\tau}(\pmb{\gamma})
	= \frac{1}{2n} \sum \limits_{i=1}^{n}( w_i^{\tau}(\pmb{\gamma})- w_i^{\tau}(\B{\beta}) )y_i \mathbf{x}_i.
	\end{equation}
	
	For $\mathbf{a}, \mathbf{b}\in \mathbb{R}^{n}$ we define the vector $\mathbf{a}*\mathbf{b} = (a_i b_i)_{i=1}^n$. Then we can rewrite Equation \eqref{diff-gradient} as
	\begin{equation} \label{diff-gradient-bis}
	\nabla g^{\tau}(\B{\beta}) - \nabla g^{\tau}(\pmb{\gamma})
	= \frac{1}{2n} \mathbb{X}^T\left[  \mathbf{y}* \left( \textbf{w}^{\tau}(\pmb{\gamma})- \textbf{w}^{\tau}(\B{\beta}) \right) \right]. 
	\end{equation}
	
	The operator norm associated to the Euclidean norm of the matrix $\mathbb{X}$ is $\| \mathbb{X} \| = \max_{ \| \mathbf{z} \|_2 = 1}  \| \mathbb{X}\textbf{z} \|_2$.
	
	Let us recall that  $\| \mathbb{X} \|^2 = \| \mathbb{X}^T \|^2 =  \| \mathbb{X}^T \mathbb{X} \| = \mu_{\max}(\mathbb{X}^T \mathbb{X} )$ corresponds to the highest eigenvalue of the matrix $\mathbb{X}^T \mathbb{X}$.
	
	Consequently, Equation \eqref{diff-gradient-bis} leads to: 
	\begin{equation} \label{first-part}
	\| \nabla L^{\tau}(\B{\beta}) - \nabla L^{\tau}(\pmb{\gamma}) \|_2 \le
	\frac{1}{2n} \| \mathbb{X} \|  \left\|   \textbf{w}^{\tau}(\pmb{\gamma})- \textbf{w}^{\tau}(\B{\beta}) \right\|_2.
	\end{equation}

	In addition, the first order necessary conditions for optimality applied to $ \textbf{w}^{\tau}(\B{\beta})$ and $\textbf{w}^{\tau}(\pmb{\gamma})$ give
	\begin{equation} \label{first-order-1}
	\sum_{i=1}^n \left\{ \frac{1}{2n}(1-y_i \mathbf{x}_i^T\B{\beta} )- \frac{\tau}{n} w_i^{\tau}(\B{\beta} ) \right\} \left\{  w_i^{\tau}(\pmb{\gamma}) - w_i^{\tau}(\B{\beta})  \right\}\le 0,
	\end{equation}
	and
	\begin{equation} \label{first-order-2}
	\sum_{i=1}^n \left\{ \frac{1}{2n}(1-y_i \mathbf{x}_i^T\pmb{\gamma}) - \frac{\tau}{n} w_i^{\tau}(\pmb{\gamma}) \right\} \left\{  w_i^{\tau}(\B{\beta}) - w_i^{\tau}(\pmb{\gamma})  \right\}\le 0.
	\end{equation}

	Then by adding Equations \eqref{first-order-1} and  \eqref{first-order-2} and rearranging the terms we have:
	\begin{align*}
	&\tau \|  \textbf{w}^{\tau}(\pmb{\gamma})- \textbf{w}^{\tau}(\B{\beta})   \|_2^2 \\
	&\le \frac{1}{2} \sum_{i=1}^n  y_i \mathbf{x}_i^T ( \B{\beta} - \pmb{\gamma} ) \left(  w_i^{\tau}(\pmb{\gamma}) - w_i^{\tau}(\B{\beta})  \right) \\
	&\le  \frac{1}{2} \| \mathbb{X}  \left( \B{\beta} - \pmb{\gamma}   \right) \|_2   \|  \mathbf{w}^{\tau}(\pmb{\gamma}) - \mathbf{w}^{\tau}(\B{\beta})   \|_2 \\
	&\le \frac{1}{2} \| \mathbb{X} \|   \|  \B{\beta} - \pmb{\gamma}  \|_2     \|  \mathbf{w}^{\tau}(\pmb{\gamma}) - \mathbf{w}^{\tau}(\B{\beta})   \|_2,
	\end{align*}
	where we have used Cauchy-Schwartz inequality. We easily derive:
	\begin{equation}\label{second-part}
	\|  \textbf{w}^{\tau}(\pmb{\gamma})- \textbf{w}^{\tau}(\B{\beta})   \|_2
	\le \frac{1}{2\tau} \| \mathbb{X} \|   \|  \B{\beta} - \pmb{\gamma}  \|_2.
	\end{equation}
	
	We conclude the proof by combining Equations \eqref{first-part} and \eqref{second-part}:
	\begin{align*}
	\| \nabla L^{\tau}(\B{\beta}) - \nabla L^{\tau}(\pmb{\gamma}) \|_2
	&\le \frac{1}{4n\tau} \| \mathbb{X} \|^2   \|  \B{\beta}- \pmb{\gamma}  \|_2\\
	&= \frac{ \mu_{\max}(n^{-1}\mathbb{X}^T \mathbb{X} )  }{4 \tau} \|  \B{\beta} - \pmb{\gamma}  \|_2.
	\end{align*}
	
	\paragraph{The case of Quantile Regression:} 
	For the quantile regression loss, the same smoothing method applies. Let us simply note that:
	\begin{align*}
	\rho_{\theta} (x) &= \max\left( (\theta-1)x, \ \theta x \right)  = \frac{1}{2}((2\theta-1)x + |x|) \\
	&=  \max_{|w| \le 1} \frac{1}{2}( (2\theta-1)x + wx).
	\end{align*}
	Hence we can immediately use the same steps than for the hinge loss -- which is a particular case of the quantile regression loss -- and define the smooth quantile regression loss $g_{\theta}^{\tau}$. Its gradient is:
	\begin{equation}
	\nabla g^{\tau}_{\theta}( \B{\beta} ) = - \frac{1}{2n} \sum \limits_{i=1}^{n}(2 \theta -1+w_i^{\tau}(\B{\beta}) )y_i \mathbf{x}_i \in \mathbb{R}^{p},
	\end{equation}
	where we still have $ w^{\tau}_i = \min\left( 1, \frac{1}{2\tau} | z_i | \right) \sign(z_i)$ but now  $z_i = y_i - \mathbf{x}_i^T\B{\beta}, \ \forall i$. The Lipschitz constant of $\nabla g^{\tau}_{\theta}$ is still given by Theorem \ref{lipschitz}.
\end{Proof}

\section {Proof of Theorem \ref{cone-condition-slope}} \label{sec: cone-condition-slope}

\begin{Proof} 
	We still assume $|h_1| \ge \ldots \ge |h_p|$. Following Equation \eqref{inf-equation-sup} it holds:
	\begin{equation}\label{basic-slope}
	S(\mathbf{h}) \le \Delta(\mathbf{h}) \le \eta | \B{\beta}^* |_S - \eta | \hat{\B{\beta}} |_S.
	\end{equation}
	
	We want to upper-bound the right-hand side of Equation \eqref{basic-slope}. We define a permutation $\phi \in \mathcal{S}_p$ such that $| \B{\beta}^* |_S = \sum_{j=1}^{k^*} \lambda_j | \beta^*_{\phi(j)} | $ and $| \hat{\beta}_{\phi(k^*+1)} | \ge \ldots \ge | \hat{\beta}_{\phi(p)} |$:
	\begin{align}\label{inequality-slope}
	\begin{split}
	\frac{1}{\eta}\Delta(\mathbf{h}) 
	&\le \sum_{j=1}^{k^*} \lambda_j | \beta^*_{\phi(j)} | - \max \limits_{\psi \in \mathcal{S}_p} \sum_{j=1}^p \lambda_j |  \hat{\beta}_{\psi(j)} |\\
	&\le \sum_{j=1}^{k^*} \lambda_j \left( | \beta^*_{\phi(j)} | -  | \hat{\beta}_{\phi(j)} | \right) -
	\sum_{j=k^*+1}^p \lambda_j  | \hat{\beta}_{\phi(j)} | \\
	&= \sum_{j=1}^{k^*} \lambda_j | h_{\phi(j)} | - \sum_{j=k^*+1}^p \lambda_j  | \hat{\beta}_{\phi(j)} | \\
	&\le \sum_{j=1}^{k^*} \lambda_j | h_{\phi(j)} | - \sum_{j=k^*+1}^p \lambda_j  | h_{\phi(j)} | .
	\end{split}
	\end{align}
	Since $\lambda$ is monotonically non decreasing: $\sum_{j=1}^{k^*} \lambda_j | h_{\phi(j)} | \le \sum_{j=1}^{k^*} \lambda_j | h_j |$.
	
	\smallskip 
	
	Because $| h_{\phi(k^*+1)} | \ge \ldots \ge | h_{\phi(p)} |$: $\sum_{j=k^*+1}^p \lambda_j  | h_j | \le \sum_{j=k^*+1}^p \lambda_j  | h_{\phi(j)} |$. 
	
	\smallskip 
	
	In addition, Equation \eqref{upper-bound-SG} from Appendix \ref{sec: appendix_cone-condition} leads to, with probability at least $1 - \frac{\delta}{2}$:
	$$| S(\mathbf{h}) | \le 14 L M \sqrt{ \frac{\log(2/ \delta)}{n} } \sum_{j=1}^p \lambda_j | h_j | \le 14 L M \sqrt{ \frac{\log(6/ \delta)}{n} } \sum_{j=1}^p \lambda_j | h_j |= \frac{\eta}{\alpha}  | \mathbf{h} |_S,$$
	where $\eta$ in defined in the statement of the theorem. Thus, combining this last equation with Equation \eqref{inequality-slope}, it holds with probability at least $1 - \frac{\delta}{2}$:
	$$- \frac{1}{\alpha}  | \mathbf{h} |_S \le  \sum_{j=1}^{k^*} \lambda_j | h_{j} | - \sum_{j=k^*+1}^p \lambda_j  | h_{j} |,$$
	which is equivalent to saying that with probability at least $1- \frac{\delta}{2}$:
	\begin{equation} 
	\sum_{j=k^*+1}^p \lambda_j | h_j |   \le \frac{\alpha+1}{\alpha-1} \sum_{j=1}^{k^*} \lambda_j | h_j |,
	\end{equation}
	that is $\mathbf{h} \in \Gamma\left(k^*,  \frac{\alpha+1}{\alpha-1} \right)$.
\end{Proof}

\section {Proof of Corollary \ref{main-results-slope}} \label{sec: results-slope}

\begin{Proof} 
	We follow the same path than in the proof of Theorem \ref{main-results}. The results of Theorem \ref{hoeffding-sup} still hold for the value of $\tau$ defined as:
	$$\tau = 14 L \mu(k^*) \sqrt{\frac{\log(3) }{n} + \frac{\log\left( 4 p/k \right) }{nk}   +  \frac{ \log\left( 6/\delta\right)  }{nk}}.$$
	As a consequence, the restricted strong convexity derived in Lemma \ref{restricted-strong-convexity} can be applied. We consequently obtain with probability at least $1-\delta$:
	\begin{equation}\label{preliminary-slope}
	\begin{split}\frac{1}{4}  \kappa(k^*) \left\{ \|\mathbf{h}\|_2^2 \wedge r(k^*) \|\mathbf{h}\|_2 \right\} 
	&\le \tau \|\mathbf{h}\|_1 + \eta \sum_{j=1}^{k^*} \lambda_j | h_{j} | - \eta \sum_{j=k^*+1}^p \lambda_j  | h_{j} |\\
	&\le \tau  \| \mathbf{h}_{S_0} \|_1 + \eta \sum_{j=1}^{k^*} \lambda_j | h_{j} |  + \tau \| \mathbf{h}_{(S_0)^c} \|_1 - \eta \sum_{j=k^*+1}^p \lambda_j  | h_{j} |.
	\end{split}
	\end{equation}
	
	We want $\tau \le \eta \lambda_{p}$, that is $14 L \mu(k^*) \sqrt{\frac{\log(3) }{n} + \frac{\log\left( 4 p/k \right) }{nk}   +  \frac{ \log\left( 6/\delta\right)  }{nk}} \le 14 \alpha L M \sqrt{\frac{\log\left( 2e \right) }{n} \log(6/ \delta)}$, which is satisfied. Hence by plugging the result in Equation \eqref{preliminary-slope} we obtain, similarly to Section  \ref{sec: appendix_cone-condition}:
	\begin{equation*}
	\begin{split}
	\frac{1}{4}  \kappa(k^*) \left\{ \|\mathbf{h}\|_2^2 \wedge r(k^*) \|\mathbf{h}\|_2 \right\} 
	&\le \tau  \| \mathbf{h}_{S_0} \|_1 + \eta \sum_{j=1}^{k^*} \lambda_j | h_{j} | \\
	&\le \tau \sqrt{k^*}  \| \mathbf{h}_{S_0}  \|_2 + \eta \sqrt{ k^* \log(2pe /k^*) } \| \mathbf{h}_{S_0}  \|_2\\
	&\le 2\eta \sqrt{k^* \log(2pe /k^*)} \| \mathbf{h}_{S_0}  \|_2 \textnormal{ since } \tau \le \eta \lambda_{p} \le \eta \lambda_{k^*} \\
	&\le 28 \alpha L M \sqrt{ \frac{k^* \log(2pe /k^*)}{n}  \log(6/ \delta)} \| \mathbf{h}  \|_2.\\
	\end{split}
	\end{equation*}
	This last equation is very similar to Equation \eqref{use-tau-lambda} in the proof of Theorem \ref{main-results}. We conclude the proof identically, and obtain a similar bound in expectation by following the proof of Corollary \ref{main-corollary}.
\end{Proof}
\end{appendices}

\end{document}